arXiv:1802.08324v1 [math.NA] 22 Feb 2018

# Combining finite element and finite difference methods for isotropic elastic wave simulations in an energy-conserving manner

Longfei Gao[a], David Keyes[a]

[a]*Division of Computer, Electrical and Mathematical Sciences and Engineering, King Abdullah University of Science and Technology, Thuwal 23955-6900, Saudi Arabia.*

---

## Abstract

We consider numerical simulation of the isotropic elastic wave equations arising from seismic applications with non-trivial land topography. The more flexible finite element method is applied to the shallow region of the simulation domain to account for the topography, and combined with the more efficient finite difference method that is applied to the deep region of the simulation domain. We demonstrate that these two discretization methods, albeit starting from different formulations of the elastic wave equation, can be joined together smoothly via weakly imposed interface conditions. Discrete energy analysis is employed to derive the proper interface treatment, leading to an overall discretization that is energy-conserving. Numerical examples are presented to demonstrate the efficacy of the proposed interface treatment.

---

## 1. Introduction

Numerical simulation of wave propagation in earth media is of vital importance in seismic studies with substantial impacts on real-life applications such as resource exploration and hazard assessment, cf. [1–4]. Elastic wave equations are often used to describe the propagation of seismic waves emitted from either earthquake sources or active sources used in land-based explorations, cf. [5, 6]. Both finite difference methods (FDMs), e.g., [7, 8], and finite element methods (FEMs), e.g., [9, 10], are widely employed to simulate the elastic waves, with the former being more popular in seismic exploration community due to its efficiency and the latter being more popular in seismological studies of larger scales where the impacts of topography and other geometric structures are prominent.

In this work, we propose to combine the more flexible finite element method and the more efficient finite difference method in order to account for non-trivial topography, yet retain computational efficiency for majority of the simulation domain. Similar attempts may be found in existing literature, e.g., [11, 12]. Novelty of this work resides in the energy-conserving interface treatment. Specifically, we split the simulation domain into two parts, namely, the shallow region and the deep region, with a straight horizontal interface as illustrated in Figure 1. The

*Email addresses:* `longfei.gao@kaust.edu.sa` (Longfei Gao), `david.keyes@kaust.edu.sa` (David Keyes)

finite element method is employed for the shallow region (i.e., top layer in Figure 1) simulation while the finite difference method is employed for the deep region (i.e., bottom layer in Figure 1) simulation. In the following, we may refer to these two regions as the FEM region and the FDM region, respectively. The 2D isotropic elastic wave equation is considered throughout this work, which is posed in the second-order displacement formulation for the FEM region and the first-order velocity-stress formulation for the FDM region, as detailed in Section 2. These choices of formulations, albeit different, are natural to carry out energy analysis for the respective discretization methods.

In the FEM region, the standard Bubnov-Galerkin approach (i.e., the solution space and the test space are identical) is adopted to discretize the isotropic elastic wave equation, leading to a semi-discretized linear system with symmetric mass and stiffness matrices. These symmetric matrices enter the definitions of the discrete kinetic and potential energies in the FEM region. In the FDM region, the finite difference operators are designed following the summation-by-parts (SBP) principle, cf. [13–16], leading to a semi-discretized linear system that mimics the behavior of the elastic wave equation from the energy analysis perspective. The two simulation regions are joined together with an interface treatment derived from discrete energy analysis. Specifically, the interface conditions are imposed weakly through penalty terms, which are often referred to as the simultaneous approximation terms (SATs) in finite difference literature, cf. [17, 15]. With carefully designed penalty terms, the overall discretization is shown to be energy-conserving for a discrete energy resembling the physical energy associated with the elastic medium.

The rest of this paper is organized as follows. In Section 2, we describe the 2D isotropic elastic wave equation and its two formulations used in the FEM region and the FDM region, respectively. In Section 3, we briefly recount the discretization procedures used in the FEM region and the FDM region, respectively, and then present the proper interface treatment to join together the two simulation regions. In Section 4, we demonstrate the efficacy of the proposed interface treatment with numerical examples. In Section 5, we remark on several relevant problems and potential future extensions. We conclude in Section 6.

## 2. Problem description

In the absence of external forces, wave propagation in elastic media can be described by the *equation of motion*:

$$\rho \ddot{u}_i = \sigma_{ij,j} \tag{1}$$

and the *constitutive relation*:

$$\sigma_{ij} = C_{ijkl} \varepsilon_{kl} \; . \tag{2}$$

Equations (1) and (2) are written in index notation where the Einstein summation convention applies to the subscripts. For the 2D case considered in this work, all the indices, e.g., $i$, $j$, $k$ and $l$ in the above equations, range from 1 to 2. Moreover, the comma in the subscripts, e.g., $\sigma_{ij,j}$, denotes spatial differentiation, the dot overhead denotes temporal differentiation while the double dots in the case of $\ddot{u}_i$ denotes double differentiation in time. The index notation is convenient for equation derivation but lacks intuitive recognition. For the remainder of this work, we associate the first index, i.e., 1, with the axis variable in the horizontal direction, i.e., $x$, and associate the second index, i.e., 2, with the axis variable in the vertical direction, i.e., $y$. Equations and symbols written in indices and axis variables will be used interchangeably in the following. For instance, $\sigma_{1j,j}$ is the same as $\frac{d\sigma_{xx}}{dx} + \frac{d\sigma_{xy}}{dy}$.

In (1) and (2), density $\rho$ and fourth-order stiffness tensor $C_{ijkl}$ are given physical parameters that characterize the elastic media. The equation of motion relates the displacement vector $u_i$ with the second-order stress tensor $\sigma_{ij}$ while the constitutive relation relates the stress tensor with the second-order strain tensor $\varepsilon_{kl}$, which can be expressed as $\varepsilon_{kl} = \frac{1}{2}(u_{k,l}+u_{l,k})$ using the displacement vector. Both stress and strain tensors are symmetric, i.e., $\sigma_{ij} = \sigma_{ji}$ and $\varepsilon_{kl} = \varepsilon_{lk}$. Symmetry in these tensors is useful in subsequent derivations. For instance, we have $\sigma_{ij}u_{j,i} = \sigma_{ji}u_{j,i} = \sigma_{ij}u_{i,j}$ given the symmetry of $\sigma_{ij}$ and consequently, $\sigma_{ij}\varepsilon_{ij} = \sigma_{ij}u_{i,j}$ . The fourth-order stiffness tensor possesses minor and major symmetries, leading to the following relations: $C_{ijkl} = C_{ijlk} = C_{jikl} = C_{klij}$. Alternatively, the constitutive relation (2) can be written in terms of the compliance tensor $S_{ijkl}$ as:

$$S_{ijkl}\sigma_{kl} = \varepsilon_{ij} \, , \tag{3}$$

where the fourth-order tensor $S_{ijkl}$ is the inverse of $C_{ijkl}$ and also possesses the same symmetries as $C_{ijkl}$ does. We will switch between the two forms of constitutive relations depending on the circumstances. In the discussion of finite difference discretization (cf. Section 3.2), the form of (3) is more convenient for energy analysis while the form of (2) is more suitable for implementation.

The kinetic energy density function associated with the elastic wave equation is defined as

$$\varrho_k = \frac{1}{2}\rho v_i v_i \, , \tag{4}$$

where $v_i = \dot{u}_i$ is the velocity vector. The potential (strain) energy density function associated with the elastic wave equation is defined as

$$\varrho_p = \frac{1}{2}\sigma_{ij}\varepsilon_{ij} \, . \tag{5}$$

Substituting (2) or (3) into (5), the potential energy density function can also be written as $\varrho_p = \frac{1}{2}\varepsilon_{ij}C_{ijkl}\varepsilon_{kl}$ or $\varrho_p = \frac{1}{2}\sigma_{ij}S_{ijkl}\sigma_{kl}$ . Integration of $\varrho_k$ and $\varrho_p$ over a given region gives the kinetic and potential energies associated with the elastic wave equation on that region, respectively. We use $e_k$ and $e_p$ to denote these two energies and append superscript $^E$ or $^D$ to indicate the region under discussion, with $^E$ for the FEM region and $^D$ for the FDM region.

In the isotropic case, the constitutive relation reduces to

$$\sigma_{ij} = \lambda\delta_{ij}\varepsilon_{kk} + 2\mu\varepsilon_{ij} \, , \tag{6}$$

parametrized by only two free parameters $\lambda$ and $\mu$, i.e., the Lamé parameters. The symbol $\delta_{ij}$ in (6) stands for the Kronecker delta, i.e., $\delta_{ij}$ equals to 1 if $i = j$; 0 otherwise.

The above brief account of elastodynamics is not meant to be thorough, but merely to serve the purpose of establishing the notation used throughout this work. For more information on the theory of elasticity, interested readers may consult [18–20]. In the following, we give formulations of the elastic wave equation in the two discretization regions.

In the FEM region, the 2D isotropic elastic wave equation is posed in the second-order displacement formulation as follows:

$$\begin{cases} \rho\ddot{u}_i &= \sigma_{ij,j} \, ; \\ \sigma_{ij} &= C_{ijkl}\varepsilon_{kl} = \lambda\delta_{ij}u_{k,k} + \mu\left(u_{i,j} + u_{j,i}\right) \, . \end{cases} \tag{7}$$

In this formulation, the displacements $u_i$, $i = 1 \dots 2$, are the sought solution variables while the stress tensor components $\sigma_{ij}$ merely serve as intermediate variables which do not necessarily appear in the actual computation. However, they will frequently appear in the upcoming derivations to simplify notations and to provide physical intuitions. We refer to (7) as the second-order displacement formulation because the temporal derivatives involved therein are second order and the displacements are the sought solution variables. We may simply refer to it as the second-order formulation if there is no ambiguity.

In the FDM region, the 2D isotropic elastic wave equation is posed in the first-order velocity-stress formulation as follows:

$$\begin{cases} \rho \dot{v}_i = \sigma_{ij,j} \,; \\ \dot{\sigma}_{ij} = C_{ijkl}\dot{\varepsilon}_{kl} = \lambda\delta_{ij}v_{k,k} + \mu(v_{i,j} + v_{j,i}) \,. \end{cases} \tag{8}$$

The second equation of (8) may also be written in the following equivalent form:

$$S_{ijkl}\dot{\sigma}_{kl} = \dot{\varepsilon}_{ij} = \tfrac{1}{2}(v_{i,j} + v_{j,i}) \tag{9}$$

to assist the discussion. In this formulation, the sought solution variables include components of both velocity $v$ and stress $\sigma$, which exhibit a reciprocal relationship in (8). Moreover, the temporal derivatives involved in (8) are first order. Hence, it is referred to as the first-order velocity-stress formulation. We may simply refer to it as the first-order formulation if there is no ambiguity. System (8) can be written in its equivalent form in (10) using the axis variables $x$ and $y$, which may be more familiar to finite difference modelers.

$$\begin{cases} \dfrac{\partial v_x}{\partial t} = \dfrac{1}{\rho}\dfrac{\partial \sigma_{xx}}{\partial x} + \dfrac{1}{\rho}\dfrac{\partial \sigma_{xy}}{\partial y}; \\ \dfrac{\partial v_y}{\partial t} = \dfrac{1}{\rho}\dfrac{\partial \sigma_{xy}}{\partial x} + \dfrac{1}{\rho}\dfrac{\partial \sigma_{yy}}{\partial y}; \\ \dfrac{\partial \sigma_{xx}}{\partial t} = (\lambda + 2\mu)\dfrac{\partial v_x}{\partial x} + \lambda\dfrac{\partial v_y}{\partial y}; \\ \dfrac{\partial \sigma_{xy}}{\partial t} = \mu\dfrac{\partial v_y}{\partial x} + \mu\dfrac{\partial v_x}{\partial y}; \\ \dfrac{\partial \sigma_{yy}}{\partial t} = \lambda\dfrac{\partial v_x}{\partial x} + (\lambda + 2\mu)\dfrac{\partial v_y}{\partial y}. \end{cases} \tag{10}$$

## 3. Methodology

In this section, we first describe the discretization methods used in the interiors of the FEM region and the FDM region, respectively, and then present the interface treatment that joins these two regions. To focus on the interface treatment, we start our discussion with the case of flat topography, i.e., the FEM region is also rectangular, as illustrated in Figure 1. However, it will become clear that the derived results apply naturally to the case of non-trivial topography, as explained in Remarks 1 and 2.

In Figure 1, the simulation domain is split into two discretization regions, i.e., the FEM region (top) and the FDM region (bottom), with a straight interface. In the FDM region, four subgrids are positioned in staggered fashion, with $\sigma_{xy}$, $v_y$ and $v_x$ each occupying one subgrid and the two normal stress components $\sigma_{xx}$ and $\sigma_{yy}$ sharing one subgrid. Grid spacing in the FDM region is the same as element width in the FEM region, both denoted as $\Delta x$. Moreover, on the interface, vertices of the elements match the grid points of the subgrid occupied by $\sigma_{xy}$. In the upcoming discussion, we may use symbols $\partial L$, $\partial R$, $\partial B$, and $\partial T$ to denote the left, right, bottom, and top boundaries, respectively, and use symbol $\partial I$ to denote the interface.

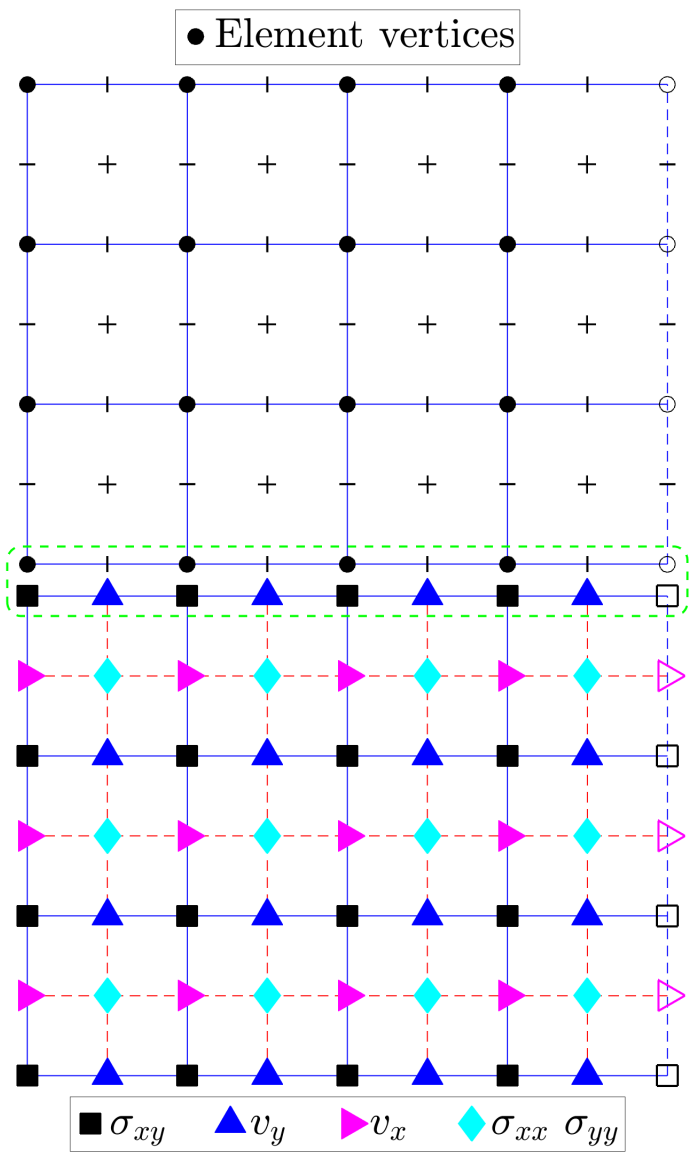

Figure 1: Geometric configuration. The simulation domain is split into the FEM region (top) and the FDM region (bottom) with a horizontal interface. The FDM region is discretized with staggered grids consisting of four subgrids. On the interface, vertices of the elements match the grid points of the $\sigma_{xy}$ subgrid.

### *3.1. Discretization in the FEM region*

Since our focus is on the interface treatment, we only give a brief account to the discretization in the FEM region here. Interested readers may consult [21–23] for more information about finite element methods.

The weak formulation of (7) can be written as: find $u_i \in \mathscr{U}_i$ such that

$$\int_\Omega \rho w_i \ddot{u}_i d_\Omega = \int_\Omega w_i \sigma_{ij,j} d_\Omega = -\int_\Omega w_{i,j} \sigma_{ij} d_\Omega + \int_{\partial\Omega} w_i \sigma_{ij} n_j d_{\partial\Omega} , \quad \forall w_i \in \mathscr{U}_i , \tag{11}$$

where $w_i$ are the test variables. The solution space and the test space are chosen to be the same in (11). Energy analysis at the continuous level can be carried out with (11) by substituting test variables $w_i$ with velocities $v_i$, leading to the following relation:

$$\int_\Omega \rho v_i \dot{v}_i d_\Omega + \int_\Omega v_{i,j} \sigma_{ij} d_\Omega = \int_{\partial\Omega} v_i \sigma_{ij} n_j d_{\partial\Omega} . \tag{12}$$

Recalling the definitions in (4) and (5), we recognize that

$$\int_\Omega \rho v_i \dot{v}_i d_\Omega = \frac{d}{dt}\left(\frac{1}{2}\int_\Omega \rho v_i v_i d_\Omega\right) = \frac{de_k^E}{dt} \tag{13}$$

and

$$\int_\Omega v_{i,j} \sigma_{ij} d_\Omega = \int_\Omega \dot{\epsilon}_{ij} \sigma_{ij} d_\Omega = \int_\Omega \sigma_{ij} S_{ijkl} \dot{\sigma}_{kl} d_\Omega = \frac{d}{dt}\left(\frac{1}{2}\int_\Omega \sigma_{ij} S_{ijkl} \sigma_{kl} d_\Omega\right) = \frac{de_p^E}{dt} , \tag{14}$$

where $e_k^E$ and $e_p^E$ are the kinetic and potential energies associated with the elastic wave equation in the FEM region, respectively. Substituting (13) and (14) into (12), we have

$$\frac{de^E}{dt} = \int_{\partial\Omega} v_i \sigma_{ij} n_j d_{\partial\Omega} \,, \tag{15}$$

where $e^E = e_k^E + e_p^E$ denotes the total energy. In other words, evolution of the total energy reduces to the behavior of the solution variables on the boundaries only.

For ease of discussion, we assume that the boundary integral in (15) vanishes on boundaries $\partial T$, $\partial L$, and $\partial R$ and reduces to the interface $\partial I$ only, using the arguments of, for instance, free surface boundary condition ($\sigma_{ij} n_j = 0$) on $\partial T$ and periodic boundary condition on $\partial L$ and $\partial R$. We note here that for the FEM region, the outward normal vector takes the values $[0, -1]^T$ on the interface $\partial I$ and therefore, (15) can be written as

$$\frac{de^E}{dt} = \int_{\partial I} v_i \sigma_{ij} n_j \, d_{\partial I} = \int_{\partial I} \left( -v_x \sigma_{xy} - v_y \sigma_{yy} \right) d_{\partial I} \,. \tag{16}$$

To discretize, we use the finite dimensional function space spanned by $\{\phi_\alpha\}_{\alpha=1}^{N^\phi}$ to approximate the space $\mathscr{U}_i$ and define the column vector $\boldsymbol{\phi}$ as $[\phi_1, \ldots, \phi_{N^\phi}]^T$. Furthermore, we use $b^1$ and $b^2$ to denote the coefficient column vectors corresponding to $\boldsymbol{\phi}$ and express the approximations of $u_1$ and $u_2$ as

$$\boldsymbol{\phi}^T b^1 \quad \text{and} \quad \boldsymbol{\phi}^T b^2 \,,$$

respectively. Finally, we sample the test variables in the order of $\{w_1, w_2\} = \{\phi_\alpha, 0\}$, $\alpha = 1, \ldots, N^\phi$ and then $\{w_1, w_2\} = \{0, \phi_\alpha\}$, $\alpha = 1, \ldots, N^\phi$. With these choices, the two area integrals in (11), i.e., $\int_\Omega \rho w_i \ddot{u}_i \, d_\Omega$ and $\int_\Omega w_{i,j} \sigma_{ij} \, d_\Omega$, reduce to the following matrix forms:

$$M\ddot{b} \quad \text{and} \quad Kb \,, \tag{17}$$

respectively, where

$$b = \begin{bmatrix} b^1 \\ b^2 \end{bmatrix}, \quad M = \begin{bmatrix} \int_\Omega \rho \boldsymbol{\phi}\boldsymbol{\phi}^T d_\Omega & \mathbf{0} \\ \mathbf{0} & \int_\Omega \rho \boldsymbol{\phi}\boldsymbol{\phi}^T d_\Omega \end{bmatrix} \tag{18}$$

and

$$K = \begin{bmatrix} \int_\Omega (\lambda + 2\mu) \boldsymbol{\phi}_{,1} \boldsymbol{\phi}_{,1}^T d_\Omega + \int_\Omega \mu \boldsymbol{\phi}_{,2} \boldsymbol{\phi}_{,2}^T d_\Omega & \int_\Omega \lambda \boldsymbol{\phi}_{,1} \boldsymbol{\phi}_{,2}^T d_\Omega + \int_\Omega \mu \boldsymbol{\phi}_{,2} \boldsymbol{\phi}_{,1}^T d_\Omega \\ \int_\Omega \mu \boldsymbol{\phi}_{,1} \boldsymbol{\phi}_{,2}^T d_\Omega + \int_\Omega \lambda \boldsymbol{\phi}_{,2} \boldsymbol{\phi}_{,1}^T d_\Omega & \int_\Omega \mu \boldsymbol{\phi}_{,1} \boldsymbol{\phi}_{,1}^T d_\Omega + \int_\Omega (\lambda + 2\mu) \boldsymbol{\phi}_{,2} \boldsymbol{\phi}_{,2}^T d_\Omega \end{bmatrix}. \tag{19}$$

In the definitions of $M$ in (18) and $K$ in (19), each area integral therein is a succinct representation of a matrix block. For instance, $\int_\Omega \rho \boldsymbol{\phi}\boldsymbol{\phi}^T d_\Omega$ shall be interpreted as

$$\begin{bmatrix} \int_\Omega \rho \cdot \phi_1 \cdot \phi_1 \, d_\Omega & \cdots & \int_\Omega \rho \cdot \phi_1 \cdot \phi_{N^\phi} \, d_\Omega \\ \vdots & \ddots & \vdots \\ \int_\Omega \rho \cdot \phi_{N^\phi} \cdot \phi_1 \, d_\Omega & \cdots & \int_\Omega \rho \cdot \phi_{N^\phi} \cdot \phi_{N^\phi} \, d_\Omega \end{bmatrix}.$$

These two matrices, i.e., $M$ and $K$, are referred to as the mass matrix and the stiffness matrix, respectively. To arrive at the form of $K$ in (19), the isotropic assumption (6) has been invoked. Discretization of the boundary integral term in (11) is left unaddressed here. As shown later in Section 3.3, it will be cancelled out by the penalty terms introduced for the interface treatment.

The discrete kinetic and potential energies in the FEM region are defined as

$$\mathcal{E}_k^E = \frac{1}{2}\dot{b}^T M \dot{b} \quad \text{and} \quad \mathcal{E}_p^E = \frac{1}{2} b^T K b \,, \tag{20}$$

respectively. Their respective correspondences to the continuous energies $e_k^E$ and $e_p^E$ can be verified straightforwardly, realizing that $\boldsymbol{\phi}^T b^1$ and $\boldsymbol{\phi}^T b^2$ are approximations to $u_1$ and $u_2$, respectively, while $\boldsymbol{\phi}^T \dot{b}^1$ and $\boldsymbol{\phi}^T \dot{b}^2$ are approximations to $v_1$ and $v_2$, respectively. Accordingly, left multiplying $\dot{b}^T$ to the two terms in (17), the resulting terms $\dot{b}^T M \ddot{b}$ and $\dot{b}^T K b$ are the discrete correspondences to $\dot{e}_k^E$ and $\dot{e}_p^E$, respectively.

**Remark 1.** *Results presented in this subsection, particularly those regarding energy analysis, extend naturally to FEM regions occupying general geometric shapes. In this case, the concept of parametric domain is often invoked to assist the implementation. Parametric domains are often of simple geometric shapes (e.g., square) so that it is easy to define finite element basis functions on them. These basis functions are then mapped to the physical domain to approximate the quantities of interest. For the above results to apply for general geometry, $\boldsymbol{\phi}$ should be understood as images of these basis functions on the FEM region (physical domain).*

*3.2. Discretization in the FDM region*

With energy analysis in mind, we consider the following form of the elastic wave equation in the FDM region:

$$\begin{cases} \rho \dot{v}_i = \sigma_{ij,j} \,; & \text{(21a)} \\ S_{ijkl}\dot{\sigma}_{kl} = \dot{\varepsilon}_{ij} = \frac{1}{2}(v_{i,j} + v_{j,i}) \,. & \text{(21b)} \end{cases}$$

Taking temporal differentiation of the total energy $e^D = e_k^D + e_p^D$ in the FDM region, we have

$$\frac{de^D}{dt} = \frac{d}{dt}\left(\frac{1}{2}\int_\Omega \left(\rho v_i v_i + \sigma_{ij} S_{ijkl}\sigma_{kl}\right) d_\Omega\right) = \int_\Omega \left(\rho v_i \dot{v}_i + \sigma_{ij} S_{ijkl}\dot{\sigma}_{kl}\right) d_\Omega \,. \tag{22}$$

Substituting the two equations of (21) into (22) and recalling the symmetry of $\sigma_{ij}$, we arrive at:

$$\frac{de^D}{dt} = \int_\Omega \left(v_i \sigma_{ij,j} + \sigma_{ij} v_{i,j}\right) d_\Omega = \int_{\partial\Omega} v_i \sigma_{ij} n_j \, d_{\partial\Omega} \,, \tag{23}$$

which is similar to the result that we obtained for the FEM region in (15), i.e., evolution of the total energy reduces to the behavior of the solution variables on the boundaries only. We aim to retain this property in the discretized system, which is often referred to as the summation-by-parts (SBP) property in the finite difference literature, cf. [13–16, 24, 25]. Design of SBP operators often omits boundary conditions, in which case the technique of simultaneous approximation terms (SATs) is often invoked to impose the boundary conditions weakly through penalty, cf. [17, 26–29].

Again, for ease of discussion, we assume that the boundary integral in (23) vanishes on boundaries $\partial B$, $\partial L$, and $\partial R$ and reduces to the interface only. We note here that for the FDM region, the outward normal vector takes the values $[0, 1]^T$ on the interface $\partial I$ and therefore, (23) can be written as

$$\frac{de^D}{dt} = \int_{\partial I} v_i \sigma_{ij} n_j d_{\partial I} = \int_{\partial I} \left( v_x \sigma_{xy} + v_y \sigma_{yy} \right) d_{\partial I} \; . \tag{24}$$

For the isotropic elastic wave equation considered in this work, we use the staggered grids demonstrated in the FDM region of Figure 1 for its discretization. This type of grid configuration is popular in seismic studies, cf. [7, 30, 31], and dates back to the Yee scheme [32]. Symbolically, we use the following system to denote the finite difference discretization of (21):

$$\begin{cases} \mathcal{A}^{V_i} \boldsymbol{\rho}^{V_i} \dot{V}_i = \mathcal{A}^{V_i} \mathcal{D}_j^{\Sigma_{ij}} \Sigma_{ij} \; ; & \text{(25a)} \\ \mathcal{A}^{\Sigma_{ij}} \boldsymbol{S}_{ijkl}^{\Sigma_{kl}} \dot{\Sigma}_{kl} = \frac{1}{2} \mathcal{A}^{\Sigma_{ij}} ( \mathcal{D}_j^{V_i} V_i + \mathcal{D}_i^{V_j} V_j ) \; . & \text{(25b)} \end{cases}$$

In (25), the Einstein summation convention applies only to indices in the subscripts. For instance, in $\mathcal{A}^{V_i} D_{\mathbf{j}}^{\Sigma_{ij}} \Sigma_{i\mathbf{j}}$, the summation only applies to the index $\mathbf{j}$ highlighted in bold font. Superscripts are appended to indicate which subgrids (or variables) these matrices are associated with. For instance, $\mathcal{A}^{V_1}$ and $\boldsymbol{\rho}^{V_1}$ are associated with the $v_x$ subgrid (or variable $v_x$). A final note on notation is that for (25b) to make sense, given any fixed $i$ and $j$, all stress components $\sigma_{kl}$ that correspond to non-zero compliance tensor components $S_{ijkl}$ need to be on the same subgrid as $\sigma_{ij}$. This can be verified for isotropic elastic wave equation and the grid configuration illustrated in the FDM region of Figure 1, as explained in Appendix A.

Solution variables $v_i$ and $\sigma_{kl}$ in (21) are approximated by column vectors $V_i$ and $\Sigma_{kl}$ in (25), respectively. These column vectors are mapped from the respective subgrids in column-wise manner. Spatial derivatives $\sigma_{ij,j}$, $v_{i,j}$ and $v_{j,i}$ are approximated by the finite difference operators $\mathcal{D}_j^{\Sigma_{ij}}$, $\mathcal{D}_j^{V_i}$ and $\mathcal{D}_i^{V_j}$, respectively.

Matrices $\mathcal{A}^{V_i}$ and $\mathcal{A}^{\Sigma_{ij}}$ may seem redundant at first glance, but will become useful in the upcoming discrete energy analysis. These matrices are referred to as the norm matrices in the SBP literature. In this work, we limit ourselves to the case of *diagonal* norm matrices. Loosely speaking, diagonal entries of $\mathcal{A}^{V_i}$ and $\mathcal{A}^{\Sigma_{ij}}$ resemble the areas that their corresponding grid points occupy, and are always positive. It is pointed out in [33] that diagonal entries of the norm matrices and their corresponding grid points provide quadrature rules for the underlying discretization domain, acting as the quadrature weights and quadrature points, respectively. Therefore, the appearances of $\mathcal{A}^{V_i}$ and $\mathcal{A}^{\Sigma_{ij}}$ in (25) can be understood as stemming from the integral operator $\int_\Omega$ that appears in the definitions of continuous energies.

Finally, matrices $\boldsymbol{\rho}^{V_i}$ and $\boldsymbol{S}_{ijkl}^{\Sigma_{kl}}$ in (25) are also *diagonal*, whose diagonal entries contain the respective discrete coefficients on subgrids indicated by their superscripts. These matrices are referred to as the coefficient matrices in the following. We remark here that diagonal matrices (e.g., norm matrices and coefficient matrices) of the same sizes commute under multiplication.

System (25) can be written in its equivalent form in (26) using the axis variables $x$ and $y$:

$$\begin{cases} \mathcal{A}^{V_x} \boldsymbol{\rho}^{V_x} \dot{V}_x = \mathcal{A}^{V_x} \left( \mathcal{D}_x^{\Sigma_{xx}} \Sigma_{xx} + \mathcal{D}_y^{\Sigma_{xy}} \Sigma_{xy} \right) ; & \text{(26a)} \\ \mathcal{A}^{V_y} \boldsymbol{\rho}^{V_y} \dot{V}_y = \mathcal{A}^{V_y} \left( \mathcal{D}_x^{\Sigma_{xy}} \Sigma_{xy} + \mathcal{D}_y^{\Sigma_{yy}} \Sigma_{yy} \right) ; & \text{(26b)} \\ \mathcal{A}^{\Sigma_{xx}} \boldsymbol{S}_{xxkl}^{\Sigma_{kl}} \dot{\Sigma}_{kl} = \mathcal{A}^{\Sigma_{xx}} \mathcal{D}_x^{V_x} V_x ; & \text{(26c)} \\ \mathcal{A}^{\Sigma_{xy}} \boldsymbol{S}_{xykl}^{\Sigma_{kl}} \dot{\Sigma}_{kl} = \frac{1}{2} \mathcal{A}^{\Sigma_{xy}} \left( \mathcal{D}_y^{V_x} V_x + \mathcal{D}_x^{V_y} V_y \right) ; & \text{(26d)} \\ \mathcal{A}^{\Sigma_{yy}} \boldsymbol{S}_{yykl}^{\Sigma_{kl}} \dot{\Sigma}_{kl} = \mathcal{A}^{\Sigma_{yy}} \mathcal{D}_y^{V_y} V_y . & \text{(26e)} \end{cases}$$

These two forms, as well as the notations therein, will be used interchangeably in the following.

Design of SBP operators has been thoroughly discussed in existing literature. Interested readers may consult [14–16] for general information and [34–36] for their design on staggered grids. In particular, 1D SBP operators on staggered grids have been devised in [36], where the concept of projection operator is introduced to deal with the situation when subgrids do not align with the boundaries. Authors of [36] also demonstrated how to construct 2D SBP operators on staggered grids using these 1D SBP operators as building blocks. Here, we omit the derivation detail and simply use the 1D SBP operators from [36] and then demonstrate how they can be applied to the simulation of isotropic elastic wave equation on the grid configuration illustrated in the FDM region of Figure 1. These 1D SBP operators are included in Appendix B so that this work can be self-contained.

Specifically, the 1D building blocks include the 1D norm matrices $\mathcal{A}_x^N$, $\mathcal{A}_x^M$, $\mathcal{A}_y^N$, and $\mathcal{A}_y^M$, the 1D finite difference operators $\mathcal{D}_x^N$, $\mathcal{D}_x^M$, $\mathcal{D}_y^N$, and $\mathcal{D}_y^M$ and the 1D identity matrices $\mathcal{I}_x^N$, $\mathcal{I}_x^M$, $\mathcal{I}_y^N$, and $\mathcal{I}_y^M$. Superscript $^N$ indicates that the operator is associated with a grid whose endpoints match the boundaries. We refer to this grid as the $N$-grid in the following. Similarly, superscript $^M$ indicates that the operator is associated with a grid that is staggered with respect to the previously mentioned $N$-grid. We refer to this staggered one as the $M$-grid in the following. To give an example, on the interface depicted in Figure 1, the $N$-grid is occupied by $\sigma_{xy}$ while the $M$-grid is occupied by $v_y$.

In the $x$-direction, we limit ourselves to the case of periodic boundary conditions in this work. In the interior, the fourth-order staggered grid central difference stencil $[1/24,\, -9/8,\, 9/8,\, -1/24]/\Delta x$ is employed, cf. [30, 37]. When approaching the left and right boundaries, the stencil is wrapped around to account for the periodic boundary condition. It can be easily verified that the resulting finite difference operators $\mathcal{D}_x^N$ and $\mathcal{D}_x^M$ satisfy the relation $\mathcal{D}_x^N + \left(\mathcal{D}_x^M\right)^T = \mathbf{0}$. Moreover, the norm matrices $\mathcal{A}_x^N$ and $\mathcal{A}_x^M$ are simply chosen as the identity matrices of the appropriate sizes, scaled by the grid spacing $\Delta x$. Consequently, we have

$$\mathcal{A}_x^N \mathcal{D}_x^M + \left(\mathcal{A}_x^M \mathcal{D}_x^N\right)^T = \mathbf{0}\,. \tag{27}$$

In the $y$-direction, the finite difference operators $\mathcal{D}_y^N$ and $\mathcal{D}_y^M$ take the forms of (B.1a) and (B.1b), respectively, while the norm matrices $\mathcal{A}_y^N$ and $\mathcal{A}_y^M$ take the forms of (B.1c) and (B.1d), respectively. For $\mathcal{D}_y^N$ and $\mathcal{D}_y^M$, the fourth-order stencil $[1/24,\, -9/8,\, 9/8,\, -1/24]/\Delta x$ is still employed for the interior, but adapts to the boundaries as in (B.1a) and (B.1b) and reduces to second order in the process.[1] Moreover, these matrices have the following property:

$$\mathcal{A}_y^N \mathcal{D}_y^M + \left(\mathcal{A}_y^M \mathcal{D}_y^N\right)^T = -\,\mathcal{E}_B^N (\mathcal{P}_B^M)^T + \mathcal{E}_I^N (\mathcal{P}_I^M)^T, \tag{28}$$

where $\mathcal{E}_B^N$, $\mathcal{E}_I^N$, $\mathcal{P}_B^M$ and $\mathcal{P}_I^M$ are column vectors, whose explicit forms are displayed in (B.3). Specifically, $\mathcal{E}_B^N$ and $\mathcal{E}_I^N$ are canonical basis vectors that select the $N$-grid values on boundary $\partial B$ and interface $\partial I$, respectively, while $\mathcal{P}_B^M$ and $\mathcal{P}_I^M$ are projection operators that project the $M$-grid values to boundary $\partial B$ and interface $\partial I$, respectively. By design, $\mathcal{P}_B^M$ and $\mathcal{P}_I^M$ provide second-order accurate projection approximations, matching the order of derivative approximations near the boundaries.

[1]Regarding this order reduction near the boundaries, interested readers may consult [38–41] and the references therein for more information.

The 2D norm matrices are built as tensor products of the 1D norm matrices as follows:

$$\begin{aligned}
&\mathcal{A}^{V_x} = \mathcal{A}_x^N \otimes \mathcal{A}_y^M, \quad \mathcal{A}^{V_y} = \mathcal{A}_x^M \otimes \mathcal{A}_y^N, \\
&\mathcal{A}^{\Sigma_{xy}} = \mathcal{A}_x^N \otimes \mathcal{A}_y^N, \quad \mathcal{A}^{\Sigma_{xx}} = \mathcal{A}^{\Sigma_{yy}} = \mathcal{A}_x^M \otimes \mathcal{A}_y^M.
\end{aligned} \tag{29}$$

The 2D finite difference operators are built as tensor products of the 1D finite difference operators and the 1D identity matrices as follows:

$$\begin{aligned}
\mathcal{D}_x^{V_x} &= \mathcal{D}_x^N \otimes \mathcal{I}_y^M, \quad \mathcal{D}_y^{V_x} = \mathcal{I}_x^N \otimes \mathcal{D}_y^M, \\
\mathcal{D}_x^{V_y} &= \mathcal{D}_x^M \otimes \mathcal{I}_y^N, \quad \mathcal{D}_y^{V_y} = \mathcal{I}_x^M \otimes \mathcal{D}_y^N, \\
\mathcal{D}_x^{\Sigma_{xy}} &= \mathcal{D}_x^N \otimes \mathcal{I}_y^N, \quad \mathcal{D}_y^{\Sigma_{xy}} = \mathcal{I}_x^N \otimes \mathcal{D}_y^N, \\
\mathcal{D}_x^{\Sigma_{xx}} &= \mathcal{D}_x^M \otimes \mathcal{I}_y^M, \quad \mathcal{D}_y^{\Sigma_{yy}} = \mathcal{I}_x^M \otimes \mathcal{D}_y^M.
\end{aligned} \tag{30}$$

The discrete kinetic and potential energies in the FDM region are defined as

$$\mathscr{E}_k^D = \frac{1}{2} V_i^T \left( \mathcal{A}^{V_i} \boldsymbol{\rho}^{V_i} \right) V_i \quad \text{and} \quad \mathscr{E}_p^D = \frac{1}{2} \Sigma_{ij}^T \left( \mathcal{A}^{\Sigma_{ij}} \boldsymbol{S}_{ijkl}^{\Sigma_{kl}} \right) \Sigma_{kl} \,, \tag{31}$$

respectively. Their respective correspondences to the continuous energies

$$e_k^D = \int_\Omega \varrho v_i v_i \, d\Omega \quad \text{and} \quad e_p^D = \int_\Omega \sigma_{ij} S_{ijkl} \sigma_{kl} \, d\Omega$$

are obvious, realizing that the norm matrices $\mathcal{A}^{V_i}$ and $\mathcal{A}^{\Sigma_{ij}}$ act as quadrature weights. The total discrete energy in the FDM region is denoted by $\mathscr{E}^D$ and defined as the sum of $\mathscr{E}_k^D$ and $\mathscr{E}_p^D$. Taking the temporal differentiation of $\mathscr{E}^D$ and substituting in (25), we have

$$\frac{d\mathscr{E}^D}{dt} = V_i^T \mathcal{A}^{V_i} \mathcal{D}_j^{\Sigma_{ij}} \Sigma_{ij} + \Sigma_{ij}^T \mathcal{A}^{\Sigma_{ij}} \mathcal{D}_j^{V_i} V_i \,. \tag{32}$$

Rewriting (32) in terms of the axis variables and collecting terms, we arrive at

$$\begin{aligned}
\frac{d\mathscr{E}^D}{dt} = \; & V_x^T \left[ \mathcal{A}^{V_x} \mathcal{D}_x^{\Sigma_{xx}} + \left( \mathcal{A}^{\Sigma_{xx}} \mathcal{D}_x^{V_x} \right)^T \right] \Sigma_{xx} \; + \; V_y^T \left[ \mathcal{A}^{V_y} \mathcal{D}_x^{\Sigma_{xy}} + \left( \mathcal{A}^{\Sigma_{xy}} \mathcal{D}_x^{V_y} \right)^T \right] \Sigma_{xy} \\
+ \; & V_x^T \left[ \mathcal{A}^{V_x} \mathcal{D}_y^{\Sigma_{xy}} + \left( \mathcal{A}^{\Sigma_{xy}} \mathcal{D}_y^{V_x} \right)^T \right] \Sigma_{xy} \; + \; V_y^T \left[ \mathcal{A}^{V_y} \mathcal{D}_y^{\Sigma_{yy}} + \left( \mathcal{A}^{\Sigma_{yy}} \mathcal{D}_y^{V_y} \right)^T \right] \Sigma_{yy} \,.
\end{aligned} \tag{33}$$

Recalling the definitions of the 2D SBP operators in (29) and (30), we have the following simplifications:

$$\mathcal{A}^{V_x} \mathcal{D}_x^{\Sigma_{xx}} + \left( \mathcal{A}^{\Sigma_{xx}} \mathcal{D}_x^{V_x} \right)^T = \left[ \mathcal{A}_x^N \mathcal{D}_x^M + \left( \mathcal{A}_x^M \mathcal{D}_x^N \right)^T \right] \otimes \mathcal{A}_y^M; \tag{34a}$$

$$\mathcal{A}^{V_y} \mathcal{D}_x^{\Sigma_{xy}} + \left( \mathcal{A}^{\Sigma_{xy}} \mathcal{D}_x^{V_y} \right)^T = \left[ \mathcal{A}_x^M \mathcal{D}_x^N + \left( \mathcal{A}_x^N \mathcal{D}_x^M \right)^T \right] \otimes \mathcal{A}_y^N; \tag{34b}$$

$$\mathcal{A}^{V_x} \mathcal{D}_y^{\Sigma_{xy}} + \left( \mathcal{A}^{\Sigma_{xy}} \mathcal{D}_y^{V_x} \right)^T = \mathcal{A}_x^N \otimes \left[ \mathcal{A}_y^M \mathcal{D}_y^N + \left( \mathcal{A}_y^N \mathcal{D}_y^M \right)^T \right]; \tag{34c}$$

$$\mathcal{A}^{V_y} \mathcal{D}_y^{\Sigma_{yy}} + \left( \mathcal{A}^{\Sigma_{yy}} \mathcal{D}_y^{V_y} \right)^T = \mathcal{A}_x^M \otimes \left[ \mathcal{A}_y^N \mathcal{D}_y^M + \left( \mathcal{A}_y^M \mathcal{D}_y^N \right)^T \right] \tag{34d}$$

for the four matrix sums appearing on the right hand side of (33), respectively. Furthermore, recalling (27), the terms in (34a) and (34b) reduce to zero blocks; recalling (28), the terms in (34c) and (34d) can be written as

$$\begin{aligned}
\mathcal{A}^{V_x}\mathcal{D}_y^{\Sigma_{xy}} + \left(\mathcal{A}^{\Sigma_{xy}}\mathcal{D}_y^{V_x}\right)^T &= \mathcal{A}_x^N \otimes \left[-\mathcal{P}_B^M\left(\mathcal{E}_B^N\right)^T + \mathcal{P}_I^M\left(\mathcal{E}_I^N\right)^T\right] \\
= -\left[\mathcal{I}_x^N \otimes \mathcal{P}_B^M\right]\cdot\mathcal{A}_x^N\cdot\left[\mathcal{I}_x^N \otimes \left(\mathcal{E}_B^N\right)^T\right] &+ \left[\mathcal{I}_x^N \otimes \mathcal{P}_I^M\right]\cdot\mathcal{A}_x^N\cdot\left[\mathcal{I}_x^N \otimes \left(\mathcal{E}_I^N\right)^T\right]
\end{aligned} \tag{35a}$$

and

$$\begin{aligned}
\mathcal{A}^{V_y}\mathcal{D}_y^{\Sigma_{yy}} + \left(\mathcal{A}^{\Sigma_{yy}}\mathcal{D}_y^{V_y}\right)^T &= \mathcal{A}_x^M \otimes \left[-\mathcal{E}_B^N\left(\mathcal{P}_B^M\right)^T + \mathcal{E}_I^N\left(\mathcal{P}_I^M\right)^T\right] \\
= -\left[\mathcal{I}_x^M \otimes \mathcal{E}_B^N\right]\cdot\mathcal{A}_x^M\cdot\left[\mathcal{I}_x^M \otimes \left(\mathcal{P}_B^M\right)^T\right] &+ \left[\mathcal{I}_x^M \otimes \mathcal{E}_I^N\right]\cdot\mathcal{A}_x^M\cdot\left[\mathcal{I}_x^M \otimes \left(\mathcal{P}_I^M\right)^T\right],
\end{aligned} \tag{35b}$$

respectively. To concentrate on the interface treatment, we assume that in the above expressions, the terms related to boundary $\partial B$ disappear in the expression of $\frac{d\mathscr{E}^D}{dt}$ due to proper boundary treatment. Eventually, (33) reduces to

$$\begin{aligned}
\frac{d\mathscr{E}^D}{dt} &= V_x^T\left[\mathcal{I}_x^N \otimes \mathcal{P}_I^M\right]\cdot\mathcal{A}_x^N\cdot\left[\mathcal{I}_x^N \otimes \left(\mathcal{E}_I^N\right)^T\right]\Sigma_{xy} \\
&+ V_y^T\left[\mathcal{I}_x^M \otimes \mathcal{E}_I^N\right]\cdot\mathcal{A}_x^M\cdot\left[\mathcal{I}_x^M \otimes \left(\mathcal{P}_I^M\right)^T\right]\Sigma_{yy}\,.
\end{aligned} \tag{36}$$

In Appendix C, we explain how to append proper penalty terms to the discretized system to account for free surface boundary condition on $\partial B$ so that those terms related to $\partial B$ do cancel out in the expression of $\frac{d\mathscr{E}^D}{dt}$.

We remark here on the resemblance between (36) and (24). Taking the first term of (36) as an example, we notice that $V_x^T\left[\mathcal{I}_x^N \otimes \mathcal{P}_I^M\right]$ projects the values of $V_x$ to the $N$-grid points on interface $\partial I$ while $\left[\mathcal{I}_x^N \otimes (\mathcal{E}_I^N)^T\right]\Sigma_{xy}$ selects the values of $\Sigma_{xy}$ at the $N$-grid points on interface $\partial I$. In addition, realizing that $\mathcal{A}_x^N$ acts as a quadrature rule for $\partial I$ with the $N$-grid points being the corresponding quadrature points, it now becomes clear that the first term of (36) is the discrete correspondence of $\int_{\partial I} v_x\sigma_{xy}\, d_{\partial I}$ in (24). Similar correspondence can be established between the second term of (36) and $\int_{\partial I} v_y\sigma_{yy}\, d_{\partial I}$ in (24).

In the following, we use symbols $\chi_I^N(\cdot)$ and $\chi_I^M(\cdot)$ to denote the restrictions (by projection, selection or other means) of a solution variable at the $N$-grid and $M$-grid points on the interface, respectively. For instance, (36) may be written as

$$\frac{d\mathscr{E}^D}{dt} = \left[\chi_I^N(v_x)\right]^T\cdot\mathcal{A}_x^N\cdot\left[\chi_I^N(\sigma_{xy})\right] + \left[\chi_I^M(v_y)\right]^T\cdot\mathcal{A}_x^M\cdot\left[\chi_I^M(\sigma_{yy})\right]. \tag{37}$$

### 3.3. *Interface treatment*

On the interface $\partial I$, we seek to impose the following interface conditions:

$$v_i^E = v_i^D\ ; \tag{38a}$$

$$\sigma_{ij}^E n_j^E + \sigma_{ij}^D n_j^D = 0\,, \tag{38b}$$

where superscripts $^E$ and $^D$ are attached to solution variables to indicate to which regions they belong, with $^E$ for the FEM region and $^D$ for the FDM region, while $n_j^E$ and $n_j^D$ are the outward normal vectors on the interface for the respective regions. (38a) stems from continuity of the elastic medium, i.e., no overlap or tear; (38b) stems from Newton's third law. On the interface, $n_j^E$ and $n_j^D$ take the values of $[0,-1]^T$ and $[0,1]^T$, respectively. Therefore, (38b) reduces to

$$\sigma_{xy}^E = \sigma_{xy}^D \quad \text{and} \quad \sigma_{yy}^E = \sigma_{yy}^D \,. \tag{39}$$

In the following, we demonstrate how to incorporate these interface conditions in the FEM and FDM discretizations in an energy-conserving manner.

First, (38b) is absorbed by the FEM discretization. Specifically, the weak formulation (11) is modified as follows: find $u_i \in \mathscr{U}_i$ such that $\forall w_i \in \mathscr{U}_i$

$$\int_\Omega \rho w_i \ddot{u}_i d_\Omega = -\int_\Omega w_{i,j}\sigma_{ij} d_\Omega + \int_{\partial I} w_i \sigma_{ij} n_j d_{\partial I} + \underline{\eta^E \int_{\partial I} w_i(\sigma_{ij}n_j + \sigma_{ij}^D n_j^D) d_{\partial I}} \,, \tag{40}$$

where (38b) is incorporated through the penalty term (underlined). Superscript $^E$ is dropped from the FEM variables since there is no ambiguity. Setting the penalty parameter $\eta^E$ to $-1$, (40) reduces to

$$\int_\Omega \rho w_i \ddot{u}_i d_\Omega = -\int_\Omega w_{i,j}\sigma_{ij} d_\Omega - \int_{\partial I} w_i \sigma_{ij}^D n_j^D d_{\partial I} \,, \tag{41}$$

which leads to the following matrix form:

$$M\ddot{b} = -Kb - p \,, \tag{42}$$

where matrices $M$ and $K$ have been given in (18) and (19), respectively, while vector $p$ takes the following form:

$$p = \begin{bmatrix} \displaystyle\int_{\partial I} \boldsymbol{\phi}\sigma_{xy}^D \, d_{\partial I} \\ \displaystyle\int_{\partial I} \boldsymbol{\phi}\sigma_{yy}^D \, d_{\partial I} \end{bmatrix}. \tag{43}$$

Left multiplying (42) with $\dot{b}^T$ and recalling the definitions of the discrete energies in (20), we have the following result:

$$\frac{d\mathscr{E}^E}{dt} = -\int_{\partial I} \left(\boldsymbol{\phi}^T \dot{b}^1\right) \cdot \sigma_{xy}^D \, d_{\partial I} - \int_{\partial I} \left(\boldsymbol{\phi}^T \dot{b}^2\right) \cdot \sigma_{yy}^D \, d_{\partial I} \tag{44}$$

regarding the discrete energy associated with (42) in the FEM region. When it comes to implementation, the line integral $\int_{\partial I}$ in (43) is usually replaced by quadrature. We use $\mathbf{x}^Q$ to denote the quadrature points for $\partial I$ and use $\mathcal{W}_x^Q$ to denote the diagonal matrix whose diagonal components are the corresponding quadrature weights. The result in (44) can now be written as:

$$\frac{d\mathscr{E}^E}{dt} = -\left[\chi_I^Q(v_x^E)\right]^T \cdot \mathcal{W}_x^Q \cdot \left[\chi_I^Q(\sigma_{xy}^D)\right] - \left[\chi_I^Q(v_y^E)\right]^T \cdot \mathcal{W}_x^Q \cdot \left[\chi_I^Q(\sigma_{yy}^D)\right] \,, \tag{45}$$

where, for instance, $\chi_I^Q(v_x^E)$ is the restriction of $v_x^E = \boldsymbol{\phi}^T \dot{b}^1$ at the quadrature points $\mathbf{x}^Q$ on the interface.

On the other hand, (38a) is absorbed by the FDM discretization. Specifically, the finite difference discretization (26) is modified as follows:

$$\begin{cases} \mathcal{A}^{V_x}\boldsymbol{\rho}^{V_x}\dot{V}_x = \mathcal{A}^{V_x}\left(\mathcal{D}_x^{\Sigma_{xx}}\Sigma_{xx} + \mathcal{D}_y^{\Sigma_{xy}}\Sigma_{xy}\right); & (46a) \\ \mathcal{A}^{V_y}\boldsymbol{\rho}^{V_y}\dot{V}_y = \mathcal{A}^{V_y}\left(\mathcal{D}_x^{\Sigma_{xy}}\Sigma_{xy} + \mathcal{D}_y^{\Sigma_{yy}}\Sigma_{yy}\right); & (46b) \\ \mathcal{A}^{\Sigma_{xx}}\boldsymbol{S}_{xxkl}^{\Sigma_{kl}}\dot{\Sigma}_{kl} = \mathcal{A}^{\Sigma_{xx}}\mathcal{D}_x^{V_x}V_x; & (46c) \\ \mathcal{A}^{\Sigma_{xy}}\boldsymbol{S}_{xykl}^{\Sigma_{kl}}\dot{\Sigma}_{kl} = \frac{1}{2}\mathcal{A}^{\Sigma_{xy}}\left(\mathcal{D}_y^{V_x}V_x + \mathcal{D}_x^{V_y}V_y\right) & (46d) \\ \qquad + \underline{\frac{1}{2}\eta_{\sigma_{xy}}^{D}\left(\mathcal{I}_x^N \otimes \mathcal{E}_I^N\right)\mathcal{A}_x^N\left[\left(\mathcal{I}_x^N \otimes \left(\mathcal{P}_I^M\right)^T\right)V_x - \chi_I^N\left(v_x^E\right)\right]}; & \\ \mathcal{A}^{\Sigma_{yy}}\boldsymbol{S}_{yykl}^{\Sigma_{kl}}\dot{\Sigma}_{kl} = \mathcal{A}^{\Sigma_{yy}}\mathcal{D}_y^{V_y}V_y & (46e) \\ \qquad + \underline{\eta_{\sigma_{yy}}^{D}\left(\mathcal{I}_x^M \otimes \mathcal{P}_I^M\right)\mathcal{A}_x^M\left[\left(\mathcal{I}_x^M \otimes \left(\mathcal{E}_I^N\right)^T\right)V_y - \chi_I^M\left(v_y^E\right)\right]}, & \end{cases}$$

where the additional penalty terms (underlined) in (46d) and (46e) result from the discretization of (38a). Superscript $^D$ is dropped from the FDM variables since there is no ambiguity. Setting penalty parameters $\eta_{\sigma_{xy}}^D$ and $\eta_{\sigma_{yy}}^D$ to $-1$ and following the procedure described in Section 3.2, we arrive at the following result

$$\frac{d\mathscr{E}^D}{dt} = \left[\chi_I^N(v_x^E)\right]^T \cdot \mathcal{A}_x^N \cdot \left[\chi_I^N(\sigma_{xy}^D)\right] + \left[\chi_I^M(v_y^E)\right]^T \cdot \mathcal{A}_x^M \cdot \left[\chi_I^M(\sigma_{yy}^D)\right] \tag{47}$$

regarding the discrete energy associated with (46) in the FDM region.

We remark here that (46) is convenient for carrying out the discrete energy analysis, but cumbersome for implementation. Instead, the penalty terms in (46d) and (46e) can be absorbed by the $y$-derivative approximations $\mathcal{D}_y^{V_x}V_x$ and $\mathcal{D}_y^{V_y}V_y$, respectively, resulting in the following modified derivative approximations:

$$\widetilde{\mathcal{D}_y^{V_x}V_x} = \left[\mathcal{D}_y^{V_x} - \mathcal{I}_x^N \otimes \left(\left(\mathcal{A}_y^N\right)^{-1}\left(\mathcal{E}_I^N\left(\mathcal{P}_I^M\right)^T\right)\right)\right]V_x + \left[\mathcal{I}_x^N \otimes \left(\left(\mathcal{A}_y^N\right)^{-1}\mathcal{E}_I^N\right)\right]\chi_I^N(v_x^E)\,; \tag{48a}$$

$$\widetilde{\mathcal{D}_y^{V_y}V_y} = \left[\mathcal{D}_y^{V_y} - \mathcal{I}_x^M \otimes \left(\left(\mathcal{A}_y^M\right)^{-1}\left(\mathcal{P}_I^M\left(\mathcal{E}_I^N\right)^T\right)\right)\right]V_y + \left[\mathcal{I}_x^M \otimes \left(\left(\mathcal{A}_y^M\right)^{-1}\mathcal{P}_I^M\right)\right]\chi_I^M(v_y^E)\,, \tag{48b}$$

respectively. With these definitions, (46) can be inverted to the more familiar form:

$$\begin{cases} \dot{V}_x = \left(\boldsymbol{\rho}^{V_x}\right)^{-1}\left(\mathcal{D}_x^{\Sigma_{xx}}\Sigma_{xx} + \mathcal{D}_y^{\Sigma_{xy}}\Sigma_{xy}\right); \\ \dot{V}_y = \left(\boldsymbol{\rho}^{V_y}\right)^{-1}\left(\mathcal{D}_x^{\Sigma_{xy}}\Sigma_{xy} + \mathcal{D}_y^{\Sigma_{yy}}\Sigma_{yy}\right); \\ \dot{\Sigma}_{xx} = (\boldsymbol{\lambda}^{\Sigma_{xx}} + 2\boldsymbol{\mu}^{\Sigma_{xx}})\mathcal{D}_x^{V_x}V_x + \boldsymbol{\lambda}^{\Sigma_{xx}}\widetilde{\mathcal{D}_y^{V_y}V_y}; \\ \dot{\Sigma}_{xy} = \boldsymbol{\mu}^{\Sigma_{xy}}\mathcal{D}_x^{V_y}V_y + \boldsymbol{\mu}^{\Sigma_{xy}}\widetilde{\mathcal{D}_y^{V_x}V_x}; \\ \dot{\Sigma}_{yy} = \boldsymbol{\lambda}^{\Sigma_{yy}}\mathcal{D}_x^{V_x}V_x + (\boldsymbol{\lambda}^{\Sigma_{yy}} + 2\boldsymbol{\mu}^{\Sigma_{yy}})\widetilde{\mathcal{D}_y^{V_y}V_y}, \end{cases} \tag{49}$$

which is more suitable for implementation. Comparing with the case without penalty terms, cf. (10), the only difference in implementing (49) is that the two $y$-derivative approximations, i.e., $\mathcal{D}_y^{V_x}V_x$ and $\mathcal{D}_y^{V_y}V_y$, need to be modified according to (48) before being used to update the stress components.

We observe from (42) and (46) that $\chi_I^Q(\sigma_{xy}^D)$ and $\chi_I^Q(\sigma_{yy}^D)$ act as external input terms for the FEM region discretization while $\chi_I^N(v_x^E)$ and $\chi_I^M(v_y^E)$ act as external input terms for the FDM region discretization. Next, we demonstrate how to construct these *external input terms* so that the overall discretization is energy-conserving in the sense that

$$\frac{d\mathcal{E}^E}{dt} + \frac{d\mathcal{E}^D}{dt} = 0\,, \tag{50}$$

i.e., the remaining terms in (45) and (47) cancel out. Given the quadrature weights $\mathcal{W}_x^Q$ and the norm matrices $\mathcal{A}_x^N$ and $\mathcal{A}_x^M$, this task reduces to properly interpolating solution values between the quadrature points $\mathbf{x}^Q$ and the $N$-grid points and $M$-grid points on the interface. The $N$-grid points and $M$-grid points are denoted by $\mathbf{x}^N$ and $\mathbf{x}^M$ hereafter, respectively.

We use symbols $\mathcal{T}_{EQ}^{DN}$ and $\mathcal{T}_{EQ}^{DM}$ to denote the interpolation operators that map values from $\mathbf{x}^Q$ to $\mathbf{x}^N$ and $\mathbf{x}^M$, respectively, and symbols $\mathcal{T}_{DN}^{EQ}$ and $\mathcal{T}_{DM}^{EQ}$ to denote the interpolation operators that map values from $\mathbf{x}^N$ and $\mathbf{x}^M$, respectively, to $\mathbf{x}^Q$. Considering the first terms in (45) and (47), they become

$$-\left[\chi_I^Q(v_x^E)\right]^T \cdot \mathcal{W}_x^Q\,\mathcal{T}_{DN}^{EQ} \cdot \left[\chi_I^N(\sigma_{xy}^D)\right] \quad \text{and} \quad \left[\chi_I^Q(v_x^E)\right]^T \cdot \left(\mathcal{T}_{EQ}^{DN}\right)^T \mathcal{A}_x^N \cdot \left[\chi_I^N(\sigma_{xy}^D)\right],$$

respectively. Similarly, the second terms in (45) and (47) become

$$-\left[\chi_I^Q(v_y^E)\right]^T \cdot \mathcal{W}_x^Q\,\mathcal{T}_{DM}^{EQ} \cdot \left[\chi_I^M(\sigma_{yy}^D)\right] \quad \text{and} \quad \left[\chi_I^Q(v_y^E)\right]^T \cdot \left(\mathcal{T}_{EQ}^{DM}\right)^T \mathcal{A}_x^M \cdot \left[\chi_I^M(\sigma_{yy}^D)\right],$$

respectively. For (50) to hold regardless of the solution state, these interpolation operators need to satisfy the following relations:

$$\mathcal{W}_x^Q\,\mathcal{T}_{DN}^{EQ} = \left(\mathcal{T}_{EQ}^{DN}\right)^T \mathcal{A}_x^N \quad \text{and} \quad \mathcal{W}_x^Q\,\mathcal{T}_{DM}^{EQ} = \left(\mathcal{T}_{EQ}^{DM}\right)^T \mathcal{A}_x^M\,. \tag{51}$$

In Appendix D, we give pairs of interpolation operators that satisfy (51) for two different quadrature rules. By design, these operators provide at least second-order accurate interpolation results, matching the order of projection operators and the order of derivative approximations near the boundaries. Design of these operators is assisted by the symbolic computing software Maple.

**Remark 2.** *Based on the above derivation, it is clear that the interface treatment presented in this subsection depends only on terms that are restricted to the interface, cf.* (45) *and* (47). *Therefore, it applies naturally to the case of non-trivial topography as long as the interface is not distorted by the mapping from parametric domain to physical domain (cf. Remark 1) beyond uniform stretching in the x-direction. A numerical example of such case is presented in Section 4.2.*

### 3.4. Full discretization

In previous subsections, we have described in detail the spatial discretizations in the FEM region and the FDM region, as well as the interface treatment, resulting in the semi-discretized systems (42) and (49). In the following, we explain how these semi-discretized systems can be numerically integrated in time. We note here that exchange of information at the interface via external input terms requires the time discretizations of these two systems to be coordinated. In this work, we demonstrate how to achieve this with the staggered leapfrog time integration scheme, which is popular in seismic studies for being time reversible and easy to implement. However, other time integration schemes can be applied in similar manners as outlined in the following.

We start with the time integration of (49) with the staggered leapfrog scheme in the FDM region. Specifically, the stress components $\sigma_{xx}$, $\sigma_{xy}$ and $\sigma_{yy}$ are discretized at integer time steps, i.e., 0, 1, 2, etc.; the velocities $v_x$ and $v_y$ are discretized at half time steps, i.e., $\frac{1}{2}$, $\frac{3}{2}$, $\frac{5}{2}$, etc. Taking $v_x$ as an example, it is updated from time step $i_t - \frac{1}{2}$ to $i_t + \frac{1}{2}$ using the following formula:

$$V_x^{(i_t+1/2)} = V_x^{(i_t-1/2)} + \Delta t \left(\boldsymbol{\rho}^{V_x}\right)^{-1} \left(\mathcal{D}_x^{\Sigma_{xx}} \Sigma_{xx}^{(i_t)} \;+\; \mathcal{D}_y^{\Sigma_{xy}} \Sigma_{xy}^{(i_t)}\right),$$

where $\Delta t$ denotes the time step length while superscripts are appended to solution vectors to indicate the time steps. Other solution variables are updated in similar manners. We note here that in order to update the stress components, the external input terms $\chi_I^N(v_x^E)$ and $\chi_I^M(v_y^E)$ need to be available at half time steps.

To have a matching time integration scheme in the FEM region, the second-order system (42) is first split as follows:

$$\begin{cases} M\dot{\xi} = -Kb - p\,; & \text{(52a)} \\ \dot{b} = \xi\,, & \text{(52b)} \end{cases}$$

with the assistance of auxiliary variable $\xi$, which, like $b$ in (18), consists of subvectors $\xi^1$ and $\xi^2$. Since $\xi = \dot{b}$, $\boldsymbol{\phi}^T\xi^1$ and $\boldsymbol{\phi}^T\xi^2$ are approximations to $v_x^E$ and $v_y^E$, respectively. This first-order system is then discretized using the same staggered leapfrog scheme as in the FDM region. Specifically, $b$ is discretized at integer time steps while $\xi$ is discretized at half time steps. The external input terms $\chi_I^Q(\sigma_{xy}^D)$ and $\chi_I^Q(\sigma_{yy}^D)$ enter (52) through vector $p$ and therefore, need to be available at integer time steps in order to update $\xi$.

The above updating procedures are sketched in the following table, where the underlined variables appear in external input terms.

| | FDM | FEM |
|---|---|---|
| $i_t \longrightarrow i_t + 1/2$ | $\left\{\sigma_{xx}^D, \sigma_{xy}^D, \sigma_{yy}^D\right\} \longrightarrow \left\{v_x^D, v_y^D\right\}$ | $\left\{b, \underline{\sigma_{xy}^D}, \underline{\sigma_{yy}^D}\right\} \longrightarrow \xi$ |
| $i_t + 1/2 \longrightarrow i_t + 1$ | $\left\{v_x^D, v_y^D, \underline{\xi}\right\} \longrightarrow \left\{\sigma_{xx}^D, \sigma_{xy}^D, \sigma_{yy}^D\right\}$ | $\xi \longrightarrow b$ |

Table 1: Updating procedures for semi-discretized systems in the FDM region and the FEM region, i.e., (49) and (42), respectively.

We observe from Table 1 that by carefully matching the updates of solution variables from both discretizations as outlined in the above, the external input terms are made available for each other at the right time instances.

## 4. Numerical examples

In the following, we corroborate the proposed interface treatment with numerical examples. For this purpose, the 2D isotropic elastic wave equation is considered, which is posed in its second-order formulation (7) for the FEM region and first-order formulation (8) for the FDM region, with the corresponding discretizations described in Sections 3.1 and 3.2, respectively, and the interface treatment presented in Section 3.3. Periodic boundary conditions are considered in the $x$-direction, while free surface boundary conditions are imposed on the top and bottom boundaries. At the beginning of the simulation, the medium is assumed to be at rest, which translates to the initial conditions that all solution components, as well as their derivatives, are zero. To drive the wave propagation, a point source is imposed on the normal stress components $\sigma_{xx}$ and $\sigma_{yy}$ (with the same temporal profile for both $\sigma_{xx}$ and $\sigma_{yy}$), mimicking an explosive source in seismic survey. The simulated vertical velocity $v_y$ at some receiver location is recorded at each time step for later comparison, which may be referred to as the seismogram.

### *4.1. Example: flat topography*

To have a better focus on the interface treatment, we first consider the case of flat topography, i.e., both the FEM region and the FDM region are rectangles, as illustrated in Figure 1. The FEM region is uniformly partitioned with 200 (horizontal) by 30 (vertical) square elements. Tensor products of quadratic Lagrange basis functions with equidistant interpolatory points are used to approximate the solution. To approximate the integrals that appear in the finite element discretization, the Gauss quadrature rule with three quadrature points on both directions per element is employed. On the other hand, the FDM region is discretized by four staggered subgrids, the outmost of which is occupied by $\sigma_{xy}$ and consists of 201 (horizontal) by 31 (vertical) grid points (including the rightmost grid column, which is voided in the simulation because of the periodic boundary conditions in the $x$-direction). The finite difference operators employed in the FDM region are described in Section 3.2. Element width in the FEM region and grid spacing in the FDM region are kept the same, both denoted as $\Delta x$. On the interface, we use the interpolation operators defined by (D.1) and (D.2) of Appendix D.

The isotropic elastic medium is parametrized by density $\rho$, compressional wave velocity $c_p$ and shear wave velocity $c_s$. The Lamé parameters in constitutive relation (6), i.e., $\lambda$ and $\mu$, are linked with these parameters via $\lambda = \rho(c_p^2 - 2c_s^2)$ and $\mu = \rho c_s^2$, respectively. In this example, we consider homogeneous medium with $\rho = 1\text{kg/m}^3$, $c_p = 2\text{m/s}$ and $c_s = 1\text{m/s}$.

The point source is placed at $5\frac{1}{2}\Delta x$ below the top boundary and $49\frac{1}{2}\Delta x$ to the right of the left boundary. Its temporal profile is chosen as the Ricker wavelet with central frequency of 5Hz and time delay of 0.25s (see Appendix E for more information). We count the maximal frequency of the source content as 12.5Hz, which corresponds to minimal wavelength of 0.08m. $\Delta x$ is chosen as 0.005m, which amounts to 16 grid points (elements) per minimal wavelength. The time step length $\Delta t$ is chosen as 5e-4s while the number of simulated time steps is chosen as 20000, which amounts to 10s in total. The $v_y$ component of the simulated solution is recorded at the receiver location, which is placed at $6\Delta x$ below the top boundary and $149\frac{1}{2}\Delta x$ to the right of the left boundary.

The recorded $v_y$ signal, referred to as the FEM-FDM result, is displayed in Figure 2. For comparison, results simulated with FEM discretization alone and FDM discretization alone are also displayed in Figure 2, referred to as the FEM result and the FDM result, respectively. We observe that all three results match extremely well with one another, which implies that with the proposed interface treatment, the mixed FEM-FDM discretization is capable of delivering accurate simulation results parallel to those produced by FEM discretization and FDM discretization alone. In particular, we do not observe any additional spurious wave packets in the FEM-FDM result that may emit from the interface, thanks to the energy-conserving principle that we adhere to when designing the interface treatment.

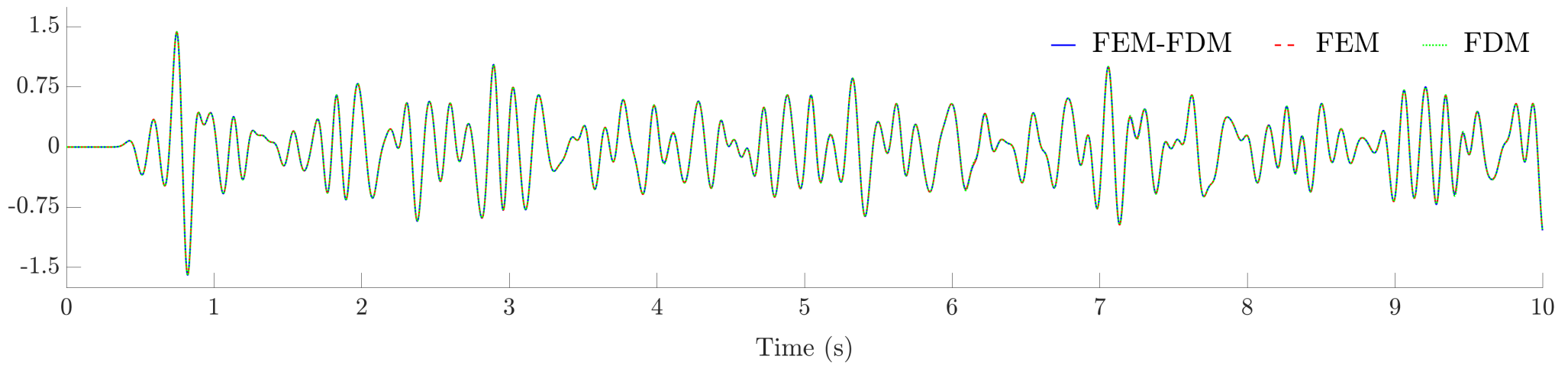

Figure 2: Seismograms simulated with mixed FEM-FDM discretization, FEM discretization alone and FDM discretization alone for the domain illustrated in Figure 1. All three superposed results match extremely well.

Next, in Figure 3, we display evolution of the total discrete energies associated with the three aforementioned simulations after the source tapers off (see Appendix E for some explanation). At each time step, the total discrete energy is calculated based on (20) for the FEM simulation and (31) for the FDM simulation. For the mixed FEM-FDM simulation, it is calculated with (20) for the FEM region and (31) for the FDM region, and then summed together. We observe that in all three cases, the total discrete energy remains constant after the source tapers off. Specifically, this confirms that the mixed FEM-FDM discretization is energy-conserving. Importance of such energy-conserving property is twofold. First, it implies stability of the simulation since the solution amplitude (measured by the discrete energy) is not allowed to grow, thus avoiding instability issues that may damage the wave simulation results (see [42] for an example). Second, the non-dissipative behavior offers additional value for seismic applications that rely on amplitude information or require long time simulations, when compared to, for instance, the Lax-Wendroff scheme, which may deliver stable simulation but yields dissipative solutions.

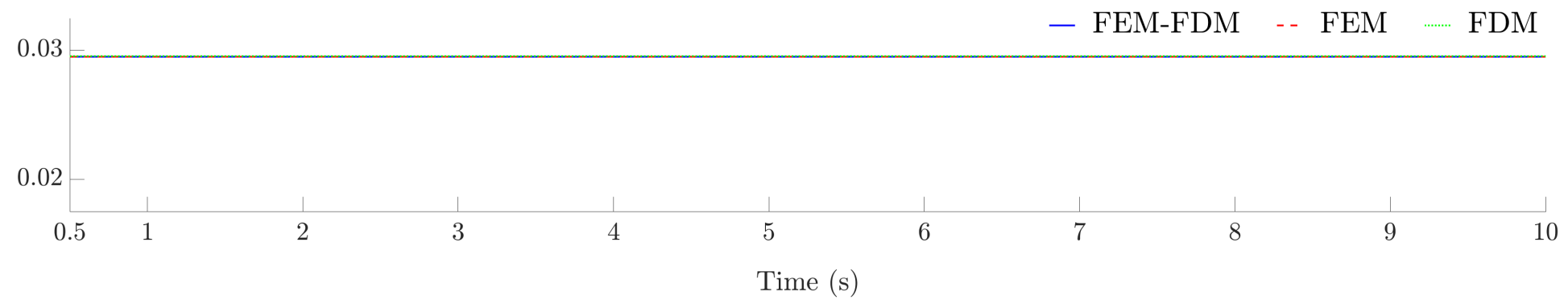

Figure 3: Evolution of the total discrete energies associated with mixed FEM-FDM simulation, FEM simulation and FDM simulation for the domain illustrated in Figure 1. In all three cases, the discrete energy remains constant after the source tapers off.

### *4.2. Example: non-trivial topography*

With this example, we demonstrate the applicability of the proposed techniques for more general cases. Specifically, we consider the simulation domain illustrated in Figure 4, where the top boundary of the FEM region (i.e., the free surface) is described by a sinusoidal function. Width and (maximal) height of the FEM region are 6250m and 500m, respectively. Amplitude of the sinusoidal function is 20% of the maximal height, leading to a variation of 200m in free surface altitude. On the other hand, the FDM region is a rectangle with width 6250m and height 1500m.

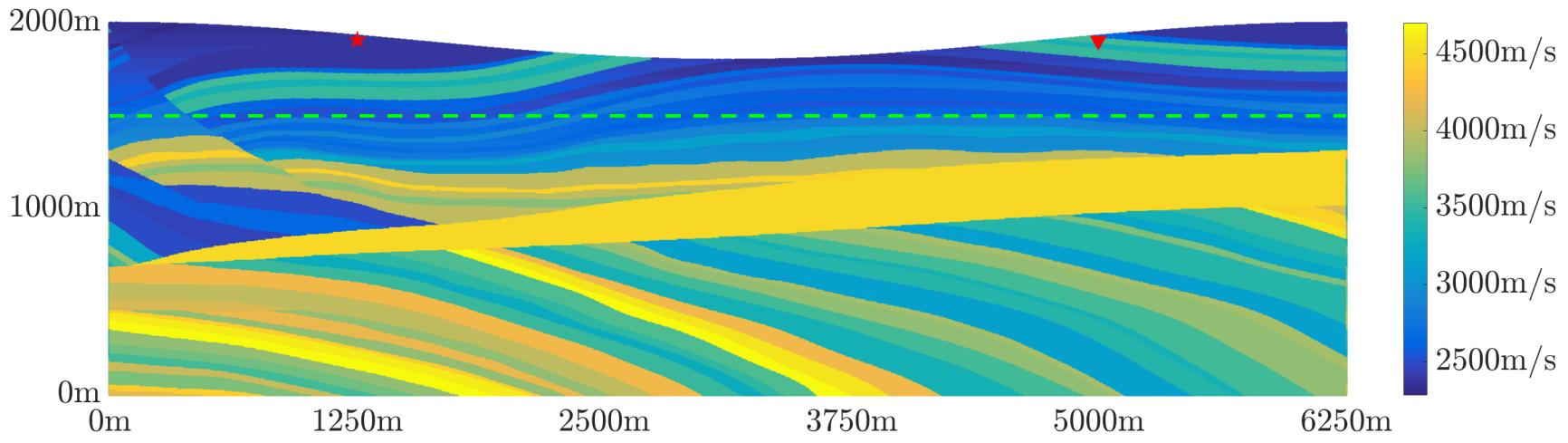

Figure 4: Illustration of the simulation domain with non-trivial topography. The FEM region (top) and the FDM region (bottom) are delineated by the green dashed line. Variation in compressional velocity $c_p$ is indicated by color in the above figure to exemplify the media complexity. Moreover, the red star and triangle indicate the source and receiver locations, respectively.

The FEM region is discretized by a mesh consisting of quadrilateral elements. We remark here that to have a convincing test, these are general irregular quadrilaterals[2], although for this simple geometry, a collection of trapezoids may suffice for the purpose of discretization. The FDM region is still discretized by the four staggered subgrids (cf. Figure 1) as in the previous example. On the interface between these two regions, we make sure that the elements (or rather, their restrictions on the interface) are uniform and that vertices of these elements match the grid points of the $\sigma_{xy}$ subgrid.

A specific variant of the finite element method, i.e., the spectral element method (see, e.g., [43, 44]) is used to assemble the finite element system for the FEM region. This approach has been popularized in global seismological studies (see, e.g., [9]). Specifically, to approximate the integrals that appear in the FEM discretization, we use the Gauss-Lobatto quadrature rule with three quadrature points on both directions per element. With Gauss-Lobatto quadrature, the quadrature points include both endpoints of the integration interval. Accordingly, bi-quadratic Lagrange basis functions whose interpolatory points coincide with these quadrature points are used to approximate the solution on each element. These choices lead to a diagonal mass matrix, which is efficient to invert and therefore, suitable for explicit dynamics. The finite difference operators used in the FDM region are the same as in the previous example. On the interface, we use the interpolation operators defined by (D.3) and (D.4) of Appendix D.

Medium parameters (i.e., $\rho$, $c_p$ and $c_s$) are cropped from the Marmousi2 model (the right bottom corner), which is a common test case for seismic studies, cf. [45]. Grid points in the FDM region are carefully aligned to match the data points in the Marmousi2 model. In the FEM region, these medium parameters are linearly interpolated to the quadrature points. To give an idea of the medium complexity, $c_p$ is indicated in Figure 4 via the variation of color. The maximum and minimum of $c_p$ are, approximately, 4700m/s and 2287m/s, respectively; the maximum and minimum of $c_s$ are, approximately, 2752m/s and 894m/s, respectively; and finally, the maximum and minimum of $\rho$ are, approximately, 2627kg/m$^3$ and 2030kg/m$^3$, respectively.

The source and receiver locations are depicted in Figure 4, which are, approximately, 23.8m and 26m below the free surface, respectively. Temporal profile of the source is the Ricker wavelet with central frequency of 5Hz and time delay of 0.25s. We count the maximal frequency of the source content as 12.5Hz, which corresponds to minimal wavelength of 71.5m, approximately. Grid spacing of the FDM region is chosen as 5m, which amounts to roughly 14.3 grid points per minimal wavelength. Element width in the FEM region, although varies from element to element, is around 5m as well. The time step length $\Delta t$ is chosen as 2e-4s while the number of simulated time steps is chosen as 30000, which amounts to 6s in total.

[2]Vertices of these quadrilaterals are included in the supplementary material in binary format (single-digits, small-endian), along with values of the medium parameters at the quadrature points.

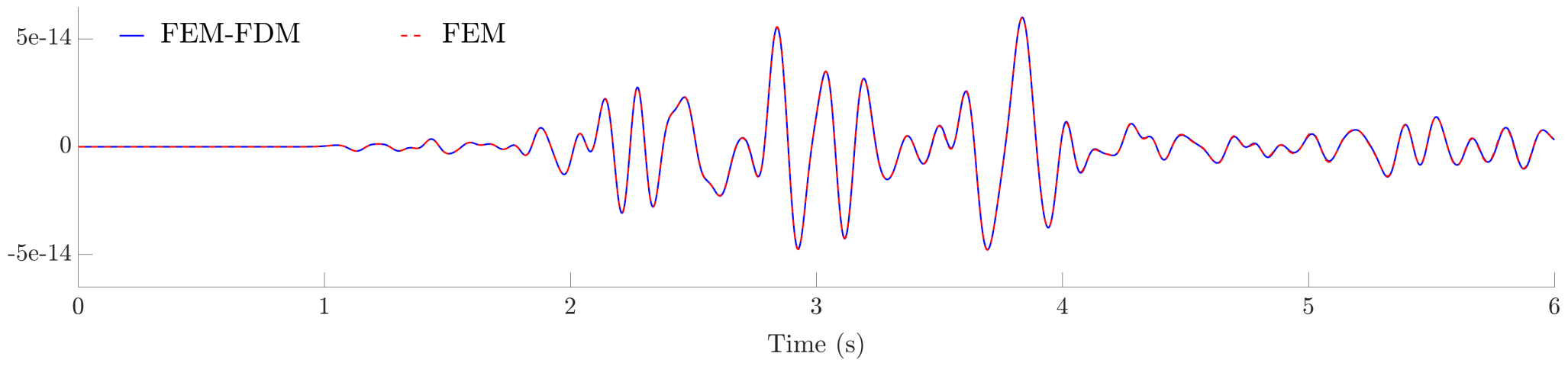

Figure 5: Seismograms simulated with mixed FEM-FDM discretization and FEM discretization alone for the domain illustrated in Figure 4.

The $v_y$ component of the simulated solution at the receiver location is recorded at every time step and displayed in Figure 5. For comparison, the result simulated with FEM discretization alone is also displayed in Figure 5. As in the previous example, we observe that the two results match very well with each other. Moreover, the total discrete energies associated with these two simulations after the source effect tapers off are displayed in Figure 6, which remain constant. Once again, this confirms that the mixed FEM-FDM discretization is energy-conserving.

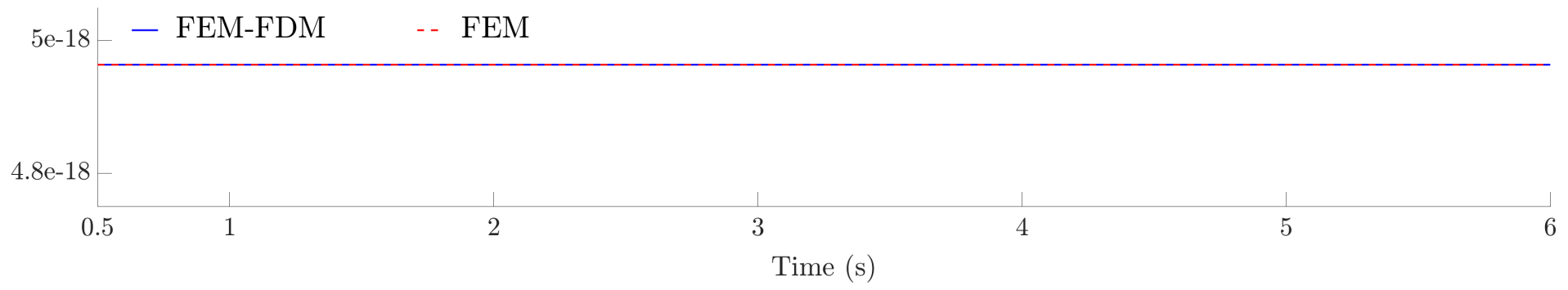

Figure 6: Evolution of the total discrete energies associated with mixed FEM-FDM simulation and FEM simulation for the domain illustrated in Figure 4.

With this example, we confirmed that the presented interface treatment can work with general meshes and heterogeneous media. This can be understood intuitively based on the derivations in Section 3. Specifically, these two factors, i.e., general meshes and heterogeneous media, are shielded from the interface treatment by 1) straight and uniformly discretized interface that makes the interface treatment immune from the mesh choices in the interior of the FEM region; 2) definitions of the discrete energies that absorb the heterogeneous media parameters such that they are absent from the interface treatment.

## 5. Discussion

In this work, we only intend to demonstrate the feasibility of combining finite element and finite difference discretizations in an energy-conserving manner for elastic wave simulations. Potential extensions and improvements remain for future exploration.

One such extension is elastic wave equations with more general constitutive relations, such as those defined on transversely isotropic media. We note here that the finite difference discretized system (25) and the discrete potential energy $\mathscr{E}_p^D$ in (31) are both written formally for general constitutive relations. Meanwhile, for the finite element discretized system (42), the constitutive relation does not appear in the penalty term $p$ associated with the interface treatment, cf. (43). Therefore, we expect the interface treatment presented in this work to transplant smoothly to more general constitutive relations.

Other branches of finite element methods can also be considered, such as the discontinuous Galerkin methods, which often start with the first-order velocity-stress formulation that is identical to the one employed in the FDM region, cf. (8) or (10). We expect this scenario to be easier to handle than the one considered in this work, where we have employed different formulations for the FEM region and the FDM region. In particular, we expect more flexibility in designing the penalty terms when the formulations employed in the two regions are the same.

Moreover, the interface treatment presented here can be naturally extended to the 3D case. Provided that restrictions of the 3D elements on the interface are uniform, the 2D interpolation operators required for information exchange on the interface can be constructed as the tensor products of their 1D counterparts, e.g., those presented in Appendix D of this work. Relations in (51) will be carried over to these 2D operators because of the properties of tensor product. In addition, since the information exchange reduces to the interface only, we expect the integration of the proposed technique to existing software, such as those from the SPECFEM project, to be marginally intrusive.

One possible improvement is in the SBP operators used in the FDM region discretization. Specifically, since we only intend to apply these SBP operators on regularly structured grids without curvilinear features, we may consider the usage of non-diagonal norm matrices to avoid the order reduction near the boundary or interface. However, existence of such operators and, if so, their incorporation in the interface treatment need to be carefully examined in future studies, which are beyond the scope of this work. Similarly, we leave the investigation of the various possible discretization orders in both regions for future studies, since the interface treatment, in particular, the design of interpolation operators, needs to be considered on a case by case basis.

## 6. Conclusion

In this work, we consider the isotropic elastic wave equations arising from land-based seismic applications. In particular, we focus on numerical techniques that enable us to combine finite element and finite difference discretizations of such equations in an energy-conserving manner. These techniques are developed following the concept of discrete energy analysis. For both finite element and finite difference discretizations, we demonstrate that with proper discretization choices, evolution of the total discrete energy, which includes a kinetic and a potential part, can be reduced to terms on the interface only. Moreover, with carefully designed interface treatment, these remaining terms on the interface can be made to cancel out each other to achieve an overall discretization that is energy-conserving. Accuracy of the proposed interface treatment and the energy-conserving property of the overall discretization are verified by numerical examples.

## 7. Acknowledgments

Gao and Keyes gratefully acknowledge the support of KAUST's Office of Sponsored Research under CCF-CAF/URF/1-2596.

## Appendix A. On the isotropic constitutive relation

The isotropic constitutive relation can be written with the axis variables as:

$$\begin{cases} \sigma_{xx} &= (\lambda + 2\mu)\varepsilon_{xx} + \lambda\varepsilon_{yy} \ ; \\ \sigma_{xy} &= 2\mu\varepsilon_{xy} \ ; \\ \sigma_{yy} &= \lambda\varepsilon_{xx} + (\lambda + 2\mu)\varepsilon_{yy} \ . \end{cases} \tag{A.1}$$

Inverting the above linear system for the strain components $\varepsilon_{xx}$, $\varepsilon_{xy}$ and $\varepsilon_{yy}$, we arrive at:

$$\begin{cases} \frac{\lambda+2\mu}{4\mu(\lambda+\mu)}\sigma_{xx} - \frac{\lambda}{4\mu(\lambda+\mu)}\sigma_{yy} &= \varepsilon_{xx} \ ; \\ \frac{1}{2\mu}\sigma_{xy} &= \varepsilon_{xy} \ ; \\ -\frac{\lambda}{4\mu(\lambda+\mu)}\sigma_{xx} + \frac{\lambda+2\mu}{4\mu(\lambda+\mu)}\sigma_{yy} &= \varepsilon_{yy} \ . \end{cases} \tag{A.2}$$

We observe from (A.2) that the non-trivial stiffness tensor components mix the two normal stress components, i.e., $\sigma_{xx}$ and $\sigma_{yy}$. For (25b) to make sense, $\sigma_{xx}$ and $\sigma_{yy}$ need to be on the same subgrid, which is the case for the grid configuration illustrated in the FDM region of Figure 1.

## Appendix B. 1D SBP operators

The following 1D SBP operators (B.1a)-(B.1d) are used as building blocks for the finite difference discretization presented in Section 3.2. They already appeared in [36] and are included here to make this work self-contained. Interested readers may consult [36] for more information.

$$\mathcal{D}_y^N = \left[\begin{array}{ccccccccccc} -79/78 & 27/26 & -1/26 & 1/78 & 0 & & & & & & \\ 2/21 & -9/7 & 9/7 & -2/21 & 0 & & & & & & \\ 1/75 & 0 & -27/25 & 83/75 & -1/25 & & & & & & \\ \hline & & 1/24 & -9/8 & 9/8 & -1/24 & & & & & \\ & & & 1/24 & -9/8 & 9/8 & -1/24 & & & & \\ & & & \ddots & \ddots & \ddots & \ddots & & & & \\ & & & & 1/24 & -9/8 & 9/8 & -1/24 & & & \\ & & & & & 1/24 & -9/8 & 9/8 & -1/24 & & \\ \hline & & & & & & 1/25 & -83/75 & 27/25 & 0 & -1/75 \\ & & & & & & 0 & 2/21 & -9/7 & 9/7 & -2/21 \\ & & & & & & 0 & -1/78 & 1/26 & -27/26 & 79/78 \end{array}\right]; \tag{B.1a}$$

$$\mathcal{D}_y^M = \left[\begin{array}{ccccccccccc} -2 & 3 & -1 & 0 & 0 & & & & & & \\ -1 & 1 & 0 & 0 & 0 & & & & & & \\ 1/24 & -9/8 & 9/8 & -1/24 & 0 & & & & & & \\ -1/71 & 6/71 & -83/71 & 81/71 & -3/71 & & & & & & \\ \hline & & 1/24 & -9/8 & 9/8 & -1/24 & & & & & \\ & & & 1/24 & -9/8 & 9/8 & -1/24 & & & & \\ & & & \ddots & \ddots & \ddots & \ddots & & & & \\ & & & & 1/24 & -9/8 & 9/8 & -1/24 & & & \\ & & & & & 1/24 & -9/8 & 9/8 & -1/24 & & \\ \hline & & & & & & 3/71 & -81/71 & 83/71 & -6/71 & 1/71 \\ & & & & & & 0 & 1/24 & -9/8 & 9/8 & -1/24 \\ & & & & & & 0 & 0 & 0 & -1 & 1 \\ & & & & & & 0 & 0 & 1 & -3 & 2 \end{array}\right]; \tag{B.1b}$$

$$\mathcal{A}_y^N = \begin{bmatrix} 7/18 & & & & & & & & & & & \\ & 9/8 & & & & & & & & & & \\ & & 1 & & & & & & & & & \\ & & & 71/72 & & & & & & & & \\ & & & & 1 & & & & & & & \\ & & & & & 1 & & & & & & \\ & & & & & & \ddots & & & & & \\ & & & & & & & 1 & & & & \\ & & & & & & & & 1 & & & \\ & & & & & & & & & 71/72 & & \\ & & & & & & & & & & 1 & \\ & & & & & & & & & & & 9/8 & \\ & & & & & & & & & & & & 7/18 \end{bmatrix}; \tag{B.1c}$$

$$\mathcal{A}_y^M = \begin{bmatrix} 13/12 & & & & & & & & & \\ & 7/8 & & & & & & & & \\ & & 25/24 & & & & & & & \\ & & & 1 & & & & & & \\ & & & & 1 & & & & & \\ & & & & & \ddots & & & & \\ & & & & & & 1 & & & \\ & & & & & & & 1 & & \\ & & & & & & & & 25/24 & & \\ & & & & & & & & & 7/8 & \\ & & & & & & & & & & 13/12 \end{bmatrix}. \tag{B.1d}$$

We note here that the matrices presented in (B.1) correspond to the case of unit grid spacing, i.e., $\Delta x = 1$. When applied to general cases, $\mathcal{D}^V$ and $\mathcal{D}^P$ need to be scaled by $1/\Delta x$ while $\mathcal{A}^P$ and $\mathcal{A}^V$ need to be scaled by $\Delta x$. With these matrices, the sum $\mathcal{A}_y^N \mathcal{D}_y^M + \left(\mathcal{A}_y^M \mathcal{D}_y^N\right)^T$, cf. (28), takes the following explicit form:

$$\begin{bmatrix} -15/8 & 5/4 & -3/8 & & & \\ & & & & & \\ & & & & & \\ & & & 3/8 & -5/4 & 15/8 \end{bmatrix}, \tag{B.2}$$

which can be written as $-\mathcal{E}_B^N(\mathcal{P}_B^M)^T + \mathcal{E}_I^N(\mathcal{P}_I^M)^T$ with $E_B^N$, $E_I^N$, $P_B^M$ and $P_I^M$ given by:

$$\mathcal{E}_B^N = \begin{bmatrix} 1 \\ 0 \\ \vdots \\ 0 \end{bmatrix}; \quad \mathcal{E}_I^N = \begin{bmatrix} 0 \\ \vdots \\ 0 \\ 1 \end{bmatrix}; \quad \mathcal{P}_B^M = \begin{bmatrix} 15/8 \\ -5/4 \\ 3/8 \\ 0 \\ \vdots \\ 0 \end{bmatrix}; \quad \mathcal{P}_I^M = \begin{bmatrix} 0 \\ \vdots \\ 0 \\ 3/8 \\ -5/4 \\ 15/8 \end{bmatrix}, \tag{B.3}$$

respectively.

## Appendix C. SATs for free surface boundary conditions

To accompany the discrete energy analysis presented throughout (32)-(35), we explain how to append the proper penalty terms to the discretized system (25) to account for free surface boundary conditions ($\sigma_{xy} = 0; \sigma_{yy} = 0$) on boundary $\partial B$ in the way that the remaining terms in $\frac{d\mathscr{E}^D}{dt}$ that are related to $\partial B$, i.e.,

$$-\,V_x^T \left[\mathcal{I}_x^N \otimes \mathcal{P}_B^M\right] \mathcal{A}_x^N \left[\mathcal{I}_x^N \otimes (\mathcal{E}_B^N)^T\right] \Sigma_{xy} \quad \text{and} \quad -\,V_y^T \left[\mathcal{I}_x^M \otimes \mathcal{E}_B^N\right] \mathcal{A}_x^M \left[\mathcal{I}_x^M \otimes (\mathcal{P}_B^M)^T\right] \Sigma_{yy}\,,$$

are cancelled out, cf. (33) and (35). Specifically, the equations used to update $V_x$ and $V_y$, i.e., (46a) and (46b), are modified as

$$\mathcal{A}^{V_x}\boldsymbol{\rho}^{V_x}\dot{V}_x = \mathcal{A}^{V_x}\left(\mathcal{D}_x^{\Sigma_{xx}}\Sigma_{xx} + \mathcal{D}_y^{\Sigma_{xy}}\Sigma_{xy}\right) + \left[\mathcal{I}_x^N \otimes \mathcal{P}_B^M\right]\mathcal{A}_x^N\left(\left[\mathcal{I}_x^N \otimes (\mathcal{E}_B^N)^T\right]\Sigma_{xy} - \mathbf{0}_x^N\right)$$

and

$$\mathcal{A}^{V_y}\boldsymbol{\rho}^{V_y}\dot{V}_y = \mathcal{A}^{V_y}\left(\mathcal{D}_x^{\Sigma_{xy}}\Sigma_{xy} + \mathcal{D}_y^{\Sigma_{yy}}\Sigma_{yy}\right) + \left[\mathcal{I}_x^M \otimes \mathcal{E}_B^N\right]\mathcal{A}_x^M\left(\left[\mathcal{I}_x^M \otimes (\mathcal{P}_B^M)^T\right]\Sigma_{yy} - \mathbf{0}_x^M\right),$$

respectively, where $\mathbf{0}_x^N$ and $\mathbf{0}_x^M$ are zero column vectors of the sizes of $\mathbf{x}^N$ and $\mathbf{x}^M$, respectively. Since $[\mathcal{I}_x^N \otimes (\mathcal{E}_B^N)^T]\Sigma_{xy}$ and $[\mathcal{I}_x^M \otimes (\mathcal{P}_B^M)^T]\Sigma_{yy}$ are restrictions of $\sigma_{xy}$ and $\sigma_{yy}$ at boundary $\partial B$, respectively, the two additional penalty terms impose relations $\sigma_{xy} = 0$ and $\sigma_{yy} = 0$, respectively.

## Appendix D. Interpolation operators

In this appendix, we give two pairs of interpolation operators that satisfy (51) for two discretization scenarios on the interface, which are considered in Sections 4.1 and 4.2, respectively.

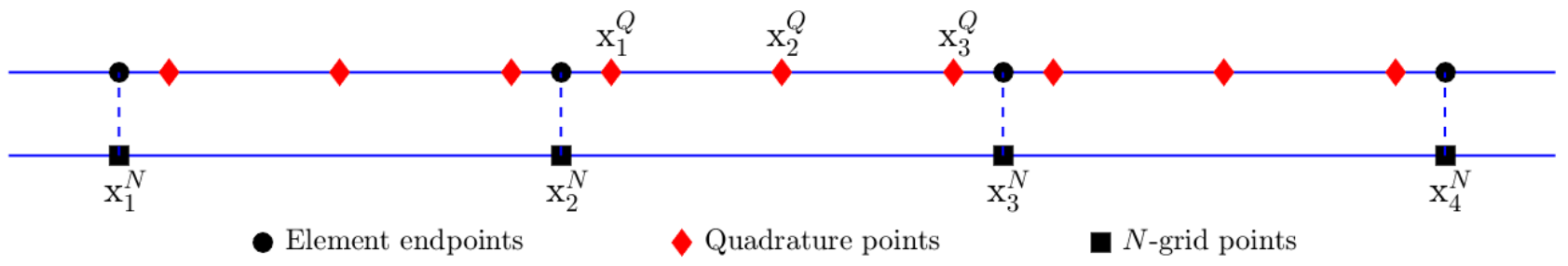

(a) Grid points involved in operator $\mathcal{T}_{DN}^{EQ}$, which interpolate from $\mathbf{x}^N$ to $\mathbf{x}^Q$.

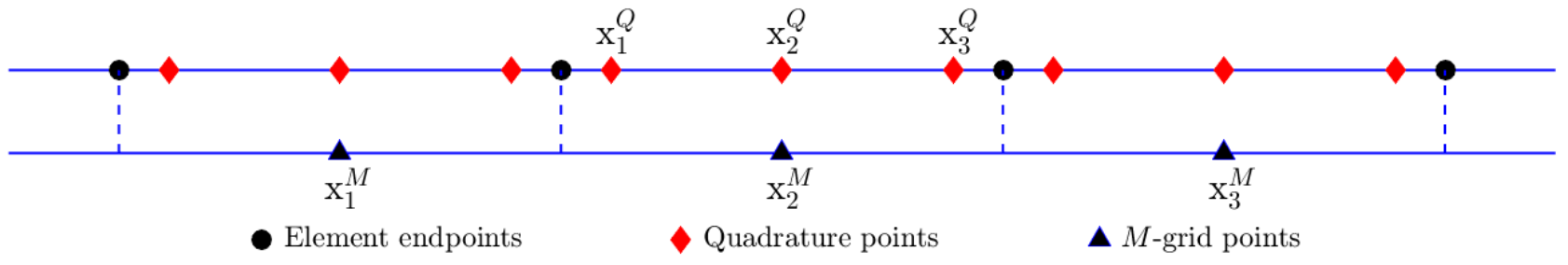

(b) Grid points involved in operator $\mathcal{T}_{DM}^{EQ}$, which interpolate from $\mathbf{x}^M$ to $\mathbf{x}^Q$.

Figure D.1: Layout of grid and quadrature points involved in the interpolations for the case of Gauss quadrature rule with three quadrature points per element. Elements are delineated by the dashed lines.

First, we consider the case where the Gauss quadrature rule with three quadrature points per element is used to approximate the line integrals in (43). Layout of grid and quadrature points involved in the interpolations is illustrated in Figure D.1. Figure D.1a corresponds to operator $\mathcal{T}_{DN}^{EQ}$, which interpolate from the $N$-grid points $\mathbf{x}^N$ to the quadrature points $\mathbf{x}^Q$. In the following, we give interpolation formulas for the three quadrature points in the middle only, since the others can be deduced based on the repeated pattern in the layout of grid points. Given the values of smooth function $f$ at $\mathrm{x}_1^N$, $\mathrm{x}_2^N$, $\mathrm{x}_3^N$ and $\mathrm{x}_4^N$, its values at $\mathrm{x}_1^Q$, $\mathrm{x}_2^Q$ and $\mathrm{x}_3^Q$ can be approximated by the following formulas:

$$\begin{aligned}
f(\mathrm{x}_1^Q) &\approx -\tfrac{1}{20}f(\mathrm{x}_1^N) + \left(\tfrac{3}{5}+\tfrac{\sqrt{15}}{10}\right)f(\mathrm{x}_2^N) + \left(\tfrac{9}{20}-\tfrac{\sqrt{15}}{10}\right)f(\mathrm{x}_3^N);\\
f(\mathrm{x}_2^Q) &\approx -\tfrac{1}{16}f(\mathrm{x}_1^N) + \tfrac{9}{16}f(\mathrm{x}_2^N) + \tfrac{9}{16}f(\mathrm{x}_3^N) - \tfrac{1}{16}f(\mathrm{x}_4^N);\\
f(\mathrm{x}_3^Q) &\approx \left(\tfrac{9}{20}-\tfrac{\sqrt{15}}{10}\right)f(\mathrm{x}_2^N) + \left(\tfrac{3}{5}+\tfrac{\sqrt{15}}{10}\right)f(\mathrm{x}_3^N) - \tfrac{1}{20}f(\mathrm{x}_4^N).
\end{aligned} \tag{D.1}$$

Similarly, Figure D.1b corresponds to operator $\mathcal{T}_{DM}^{EQ}$, which interpolates from the $M$-grid points $\mathbf{x}^M$ to the quadrature points $\mathbf{x}^Q$. The following formulas are used to interpolate from $\mathrm{x}_1^M$, $\mathrm{x}_2^M$, $\mathrm{x}_3^M$ to $\mathrm{x}_1^Q$, $\mathrm{x}_2^Q$ and $\mathrm{x}_3^Q$:

$$\begin{aligned} f(\mathrm{x}_1^Q) &\approx \left(\tfrac{3}{40}+\tfrac{\sqrt{15}}{20}\right) f(\mathrm{x}_1^M) + \tfrac{17}{20} f(\mathrm{x}_2^M) + \left(\tfrac{3}{40}-\tfrac{\sqrt{15}}{20}\right) f(\mathrm{x}_3^M); \\ f(\mathrm{x}_2^Q) &\approx f(\mathrm{x}_2^M); \\ f(\mathrm{x}_3^Q) &\approx \left(\tfrac{3}{40}-\tfrac{\sqrt{15}}{20}\right) f(\mathrm{x}_1^M) + \tfrac{17}{20} f(\mathrm{x}_2^M) + \left(\tfrac{3}{40}+\tfrac{\sqrt{15}}{20}\right) f(\mathrm{x}_3^M). \end{aligned} \tag{D.2}$$

By repeating the formulas in (D.1) and (D.2) for all quadrature points in $\mathbf{x}^Q$, one obtains the operators $\mathcal{T}_{DN}^{EQ}$ and $\mathcal{T}_{DM}^{EQ}$, respectively. Their counterparts, i.e., $\mathcal{T}_{EQ}^{DN}$ and $\mathcal{T}_{EQ}^{DM}$, can be deduced using the relations in (51).

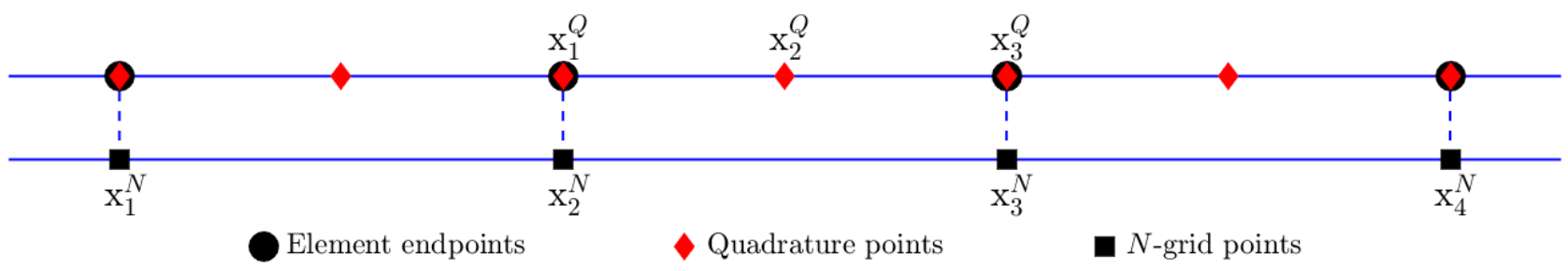

(a) Grid points involved in operator $\mathcal{T}_{DN}^{EQ}$, which interpolate from $\mathbf{x}^N$ to $\mathbf{x}^Q$.

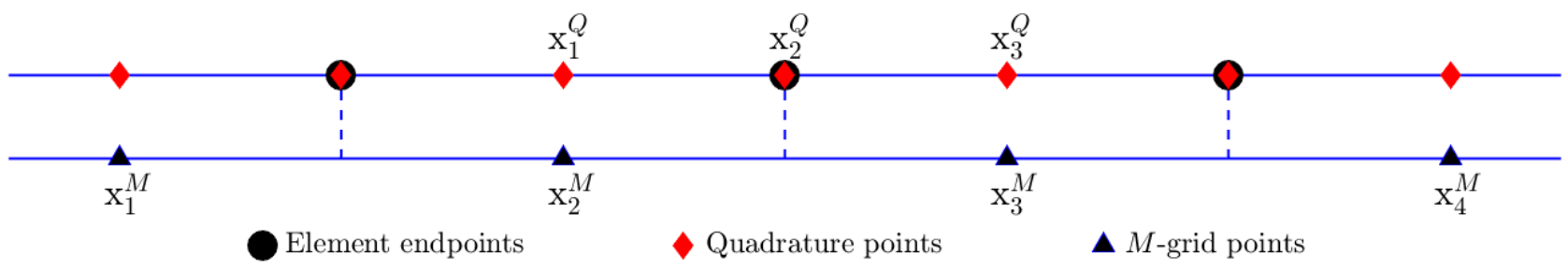

(b) Grid points involved in operator $\mathcal{T}_{DM}^{EQ}$, which interpolate from $\mathbf{x}^M$ to $\mathbf{x}^Q$.

Figure D.2: Layout of grid and quadrature points involved in the interpolations for the case of Gauss-Lobatto quadrature rule with three quadrature points per element. Elements are delineated by the dashed lines. With Gauss-Lobatto quadrature rule, element endpoints are included in the quadrature points.

Next, we consider the case where the Gauss-Lobatto quadrature rule with three quadrature points per element is used to approximate the line integrals in (43). Layout of grid and quadrature points involved in the interpolations is illustrated in Figure D.2. We note here that unlike the Gauss quadrature rule where all quadrature points are in the interior of the elements, the Gauss-Lobatto quadrature rule has the element endpoints included in the quadrature points. Figure D.2a corresponds to operator $\mathcal{T}_{DN}^{EQ}$, which interpolate from $\mathbf{x}^N$ to $\mathbf{x}^Q$ with the following formulas:

$$\begin{aligned} f(\mathrm{x}_1^Q) &\approx f(\mathrm{x}_2^N); \\ f(\mathrm{x}_2^Q) &\approx -\tfrac{1}{16} f(\mathrm{x}_1^N) + \tfrac{9}{16} f(\mathrm{x}_2^N) + \tfrac{9}{16} f(\mathrm{x}_3^N) - \tfrac{1}{16} f(\mathrm{x}_4^N); \\ f(\mathrm{x}_3^Q) &\approx f(\mathrm{x}_3^N). \end{aligned} \tag{D.3}$$

Figure D.2b corresponds to operator $\mathcal{T}_{DM}^{EQ}$, which interpolates from $\mathbf{x}^M$ to $\mathbf{x}^Q$ with the following

formulas:

$$
\begin{aligned}
f(\mathrm{x}_1^Q) &\approx f(\mathrm{x}_2^M); \\
f(\mathrm{x}_2^Q) &\approx -\tfrac{1}{16}f(\mathrm{x}_1^M) + \tfrac{9}{16}f(\mathrm{x}_2^M) + \tfrac{9}{16}f(\mathrm{x}_3^M) - \tfrac{1}{16}f(\mathrm{x}_4^M); \\
f(\mathrm{x}_3^Q) &\approx f(\mathrm{x}_3^M).
\end{aligned}
\tag{D.4}
$$

By design, all interpolation operators mentioned above, including $\mathcal{T}_{DN}^{EQ}$, $\mathcal{T}_{DM}^{EQ}$, $\mathcal{T}_{EQ}^{DN}$, and $\mathcal{T}_{EQ}^{DM}$ associated with both scenarios illustrated in Figures D.1 and D.2, provide at least second-order accurate interpolation results. These operators are derived by solving for the coefficients that satisfy the constraints demanded by accuracy and relations in (51), using the symbolic computing software Maple. Situations involving more complicated grid layouts can be handled in similar manner.

## Appendix E. Ricker wavelet (temporal profile of the source)

For the numerical examples in Section 4, a Ricker wavelet is used as the temporal profile of the point source. The standard Ricker wavelet is defined as $A(t) = (1 - 2\pi^2 f^2 t^2)e^{-\pi^2 f^2 t^2}$, which is the second derivative of a Gaussian function, where $f$ is referred to as the peak (central) frequency. Figure E.1a displays the standard Ricker wavelet for $f = 5$Hz, which is symmetric with respect to the $y$-axis. Its amplitude approaches zero as $t \to \pm\infty$ and diminishes quickly outside of a small window.

Since the simulation usually starts at $t = 0$s, it is desirable to shift the wavelet to the right so that the cutoff amplitude at $t = 0$s is negligible. This is achieved by introducing a time delay $T_0$ and variable transformation $\hat{t} = t - T_0$ so that $A(t) = (1 - 2\pi^2 f^2 \hat{t}^2)e^{-\pi^2 f^2 \hat{t}^2}$. Figure E.1b displays the shifted (delayed) Ricker wavelet for $T_0 = 0.25$s. The amount of time delay is empirical and depends on the central frequency. It needs to be larger for smaller central frequency as the shape of the wavelet is wider.

Nevertheless, once the time delay $T_0$ is determined, it is obvious from symmetry of the wavelet that after $2T_0$, the source tapers off, i.e., amplitude of the source becomes negligible. Since all analysis conducted in this work is with respect to the homogeneous wave equation (i.e., without external source terms), the energy-conserving property is only valid when external source terms are absent. Therefore, to demonstrate the energy-conserving property, the discrete energy is only displayed after the source tapers off (i.e., after $2T_0$) in Figures 3 and 6.

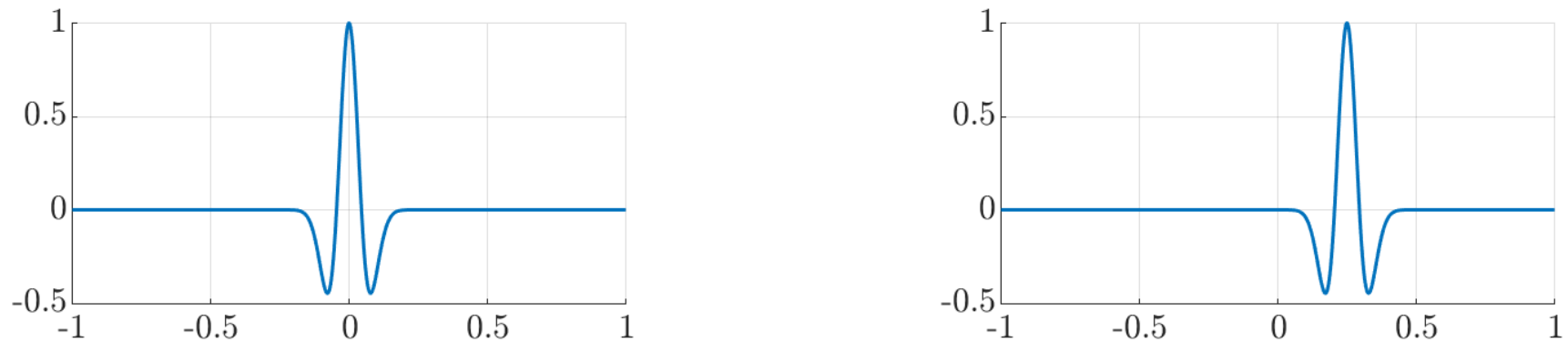

(a) Standard. (b) Shifted.

Figure E.1: Ricker wavelet. The left figure shows the standard Ricker wavelet given by formula $A(t) = (1 - 2\pi^2 f^2 t^2)e^{-\pi^2 f^2 t^2}$ with $f = 5$Hz. In the right figure, it is shifted to the right for $T_0 = 0.25$s.